\documentclass[11pt]{article}
\usepackage{amsmath, amsthm, amssymb, amsfonts}
\usepackage{graphicx}
\usepackage{subcaption}
\usepackage{setspace}
\usepackage{geometry}
\usepackage{float}
\usepackage{placeins}
\usepackage[nocompress]{cite}
\usepackage{hyperref}
\usepackage[utf8]{inputenc}
\usepackage[english]{babel}
\usepackage{framed}
\usepackage{comment}
\usepackage{etoolbox}
\usepackage{algorithm}
\usepackage{algpseudocode}
\usepackage{tikz}
\usetikzlibrary{arrows.meta,positioning,shapes.geometric}

\numberwithin{equation}{section}

\theoremstyle{plain}
\newtheorem{theorem}{Theorem}[section]
\newtheorem{proposition}[theorem]{Proposition}

\hypersetup{
	hidelinks
}

\newtheorem{corollary}[theorem]{Corollary}
\newtheorem{assumption}[theorem]{Assumption}

\theoremstyle{plain}
\newtheorem{remark}[theorem]{Remark}

\usepackage{url}

\begin{document}

\title{\bf A Validation-Driven Bilevel Optimization Method for Mesh-Free PDE Solvers and Its Error Analysis}
\author{Yongcun Song\thanks{Division of Mathematical Sciences, School of Physical and Mathematical Sciences, Nanyang Technological University, 21
Nanyang Link, 637371, Singapore. yongcun.song@ntu.edu.sg. This author was supported by the NTU Start-Up Grant.},
\and
Shangzhi Zeng\thanks{National Center for Applied Mathematics Shenzhen \& Department of Mathematics, Southern University of Science and
Technology, Shenzhen 518055, Guangdong, China. zengsz@sustech.edu.cn. This author was supported by the National Natural
Science Foundation of China (12501429) and the Shenzhen Fundamental Research Program (20250530150024003).},
\and Jin Zhang\thanks{Department of Mathematics \& National Center for Applied Mathematics Shenzhen, Southern University of Science and
Technology, Shenzhen 518055, Guangdong, China. zhangj9@sustech.edu.cn. This author was supported by the National Key R
\& D Program of China (2023YFA1011400) and National Natural Science Foundation of China (12326605, 12222106).},
\and Lvgang Zhang\thanks{Department of Mathematics, Southern University of Science and Technology, Shenzhen 518055, Guangdong, China. 12443001@mail.sustech.edu.cn.}}
\maketitle

\begin{abstract}
Mesh-free partial differential equation (PDE) solvers provide a flexible
alternative to mesh-based discretizations by constructing discrete systems
directly from residuals at sampled points, thereby avoiding mesh generation and
maintenance.
Such solvers
typically involve two related tasks: fitting coefficients in a parameterized
trial representation and selecting solver-design parameters
governing the representation, training-residual weighting, and
regularization.
These parameters strongly affect solution accuracy, yet
existing selection strategies are largely tailored to particular solver
classes. To address this limitation, we formulate a unified validation-driven bilevel optimization framework applicable across different solver classes, in which
the lower level fits coefficients using training collocation residuals, while
the upper level selects solver-design parameters using held-out validation
residuals. Conceptually, this bilevel problem can be solved by
implicit hypergradient methods.
For nonlinear PDEs, however, solving the nonlinear
lower-level problem to high accuracy is costly, while the adjoint
system is typically ill-conditioned and indefinite. 
To address this issue, we propose an efficient implicit hypergradient method that
approximates the lower-level response using finite, warm-started damped
Gauss--Newton (GN) updates. For the adjoint solve, we use the damped GN matrix
to compute an initial adjoint and to precondition a finite number of minimum
residual iterations. To
quantify the errors introduced by these approximations, we derive local
hypergradient error bounds that separate the contributions of the finite
lower-level solve, adjoint damping, and the final adjoint residual.
Combining this analysis with validation sampling bounds and residual
stability, we relate validation residuals to PDE solution accuracy and establish
error bounds for the computed PDE solution. Numerical experiments on several nonlinear PDEs, including a
high-dimensional problem, demonstrate the accuracy and 
efficiency of the proposed method across Kansa, random feature and extreme
learning machine, and physics-informed neural network solvers.
\end{abstract}

\section{Introduction}

Classical mesh-based discretizations, including finite-difference,
finite-volume, and finite-element methods, are foundational tools for solving
partial differential equations (PDEs), supported by mature theoretical frameworks and computational paradigms. Their reliance on meshes, however, creates
important challenges in high-dimensional problems and on complex or evolving
geometries. For example, the number of degrees of freedom in
conventional discretizations grows rapidly with dimension
\cite{petersdorff2004numerical,weinan2022algorithms}. Mesh generation and
maintenance are often burdensome on complex or evolving geometries
\cite{suchde2023point,ito2013challenges}. Tracking moving interfaces often requires
repeated remeshing, which adds computational cost and introduces errors when
transferring the solution between meshes \cite{li2002meshfree,gawlik2015unified}.
Residual-based mesh-free PDE solvers offer an alternative by forming discrete systems directly from residuals evaluated at sampled points, without generating or maintaining a mesh.
Such solvers have been developed for high-dimensional PDEs
\cite{sirignano2018dgm,wang2024extreme,tang2023pinns,brugiapaglia2025physics} and permit flexible placement and
adaptation of computational points on complex or evolving geometries
\cite{tominec2021unfitted,suchde2019fully,slak2019generation,adil2023ale}.
Numerical studies of several residual-based mesh-free solver classes have
reported satisfactory accuracy and favorable convergence behavior
\cite{fornberg2015solving,chen2022bridging,dong2021local}.
Building on these advantages, mesh-free solvers have been successfully
applied in areas such as geophysical modeling \cite{long2019forward}, groundwater
flow and contaminant transport \cite{seyedpour2019optimal}, image-based cardiac
electrophysiology \cite{mountris2023cardiac}, and evolving-surface modeling of
tumor growth and grain-boundary motion \cite{adil2023ale}.
\subsection{Problem Setting}
\label{sec:residual-coefficient-fitting}

For \(d\geq 1\), let \(\Omega\subset\mathbb R^d\) be a
bounded domain with boundary \(\partial\Omega\).
We consider PDEs of the form
\begin{equation*}
    \begin{cases}
        \mathcal L[u]=f, & \text{in }\Omega,\\
        \mathcal B[u]=g, & \text{on }\partial\Omega,
    \end{cases}
\end{equation*}
where \(\mathcal L\) is a possibly nonlinear differential operator,
\(\mathcal B\) is a boundary operator, and \(f\) and \(g\) are prescribed data
in suitable spaces. We assume that the
problem admits a solution \(u\) in an admissible function class on which
both operators are well defined.
For notational brevity, we write the problem in stationary form; for
time-dependent PDEs, initial conditions can be treated
analogously to boundary constraints.

To approximate \(u\), we use the finite-dimensional trial representation
\[
    \widehat{u}(\boldsymbol x;\boldsymbol\beta,\boldsymbol\omega):=
    \sum_{j=1}^{M}
    \beta_j\psi_j(\boldsymbol x;\boldsymbol\omega).
\]
Here \(\boldsymbol\omega\in\mathbb R^{n_\omega}\) denotes the solver-design
parameters governing the trial functions
\(\psi_j(\cdot;\boldsymbol\omega)\), and may also include
training-residual weights or
regularization parameters. For fixed \(\boldsymbol\omega\), the coefficient vector
\(\boldsymbol\beta\in\mathbb R^M\) determines the approximate solution within
the resulting trial space and is typically obtained by minimizing a loss
built from PDE and boundary residuals evaluated at collocation points.
Varying \(\boldsymbol\omega\) can therefore alter both the trial
representation and the coefficient-fitting problem, affecting approximation
accuracy, conditioning, and fitting stability.

In the mesh-free methods considered here, the trial functions
\(\{\psi_j(\cdot;\boldsymbol\omega)\}_{j=1}^M\) are constructed from kernels or
features independently of a mesh, unlike the mesh-based basis functions in
finite-element methods. Different constructions of
\(\psi_j(\cdot;\boldsymbol\omega)\) give rise to different solver classes, of
which the following three are representative.

\begin{itemize}
    \item \textbf{Kansa methods}, a class of radial basis function (RBF)
    collocation methods, use kernel translates as trial functions
    \cite{kansa1990multiquadrics,fasshauer1999solving,
    schaback2007convergence,ling2008stable,
    schaback2016all,bohmer2020nonlinear}.
    For prescribed centers
    \(\{\boldsymbol c_j\}_{j=1}^{M}\subset\Omega\cup\partial\Omega\) and a
    parameterized kernel \(K\), the trial functions are
    \(\psi_j(\boldsymbol x;\boldsymbol\omega)
    =K(\boldsymbol x,\boldsymbol c_j;\boldsymbol\omega)\),
    \(j=1,\ldots,M\).
    A standard example is the Gaussian kernel
    \(K(\boldsymbol x,\boldsymbol c_j;\boldsymbol\omega)
    =\exp\!\bigl(
    -\|\boldsymbol x-\boldsymbol c_j\|_2^2/
    \ell(\boldsymbol\omega)^2
    \bigr)\), where the lengthscale \(\ell(\boldsymbol\omega)>0\) parameterized
    by \(\boldsymbol\omega\) controls the shape of each kernel translate.
    In this construction, the centers specify where the kernels are placed,
    whereas the collocation points specify where the residuals are evaluated,
    and the two sets need not coincide.

    \item
    \textbf{Random feature methods} (RFM)
    \cite{rahimi2007random,nelsen2021random,mei2022generalization,
    chen2022bridging} and \textbf{extreme learning machine} (ELM)-type solvers
    \cite{huang2006extreme,dwivedi2020physics,calabro2021extreme,
    dong2021local,dong2022computing}
    generate and freeze nonlinear features before fitting their coefficients.
    For a model with one hidden layer and feature scale \(R_m>0\), the \(j\)th
    feature is
    \(\psi_j(\boldsymbol x;R_m)
    =\sigma\!\left(R_m(\widetilde{\boldsymbol w}_j\cdot\boldsymbol x
    +\widetilde b_j)\right)\),
    where \(\sigma:\mathbb R\to\mathbb R\) is a prescribed nonlinear function.
    The base variables \(\widetilde{\boldsymbol w}_j\in\mathbb R^d\) and
    \(\widetilde b_j\in\mathbb R\) are sampled once and then held fixed, so
    varying \(R_m\) changes the finite feature map. In this
    construction, the vector \(\boldsymbol\beta\) contains the fitted coefficients,
    while \(\boldsymbol\omega\) may include \(R_m\) and other selected parameters of
    the feature-generation rule, as well as training-residual weights or regularization parameters.

    \item \textbf{Physics-informed neural network (PINN)-type methods} use neural-network features as trial functions
    \cite{raissi2019physics}. When the final linear layer is separated from the
    preceding hidden layers, the network can be written as
    \(\sum_{j=1}^{M}\beta_j\psi_j(\boldsymbol x;\boldsymbol\omega)\),
    where \(\boldsymbol\beta\) contains the final-layer coefficients,
    while \(\boldsymbol\omega\) includes the hidden-layer parameters that generate
    the features and may also include training-residual weights or regularization
    parameters. For example, the Iterative PINNs with Fourier-Enhanced
    Features (IFeF-PINN) method provides such a realization
    \cite{wu2026iterative}.

\end{itemize}

\subsection{Related Work and Motivations}
\label{subsec:related-work}

The mesh-free PDE solvers described above involve two coupled tasks:
fitting the coefficients \(\boldsymbol\beta\) for given solver-design
parameters \(\boldsymbol\omega\), and selecting \(\boldsymbol\omega\) to
improve the resulting approximate solution. In this subsection, we first review
related work on Kansa, RFM/ELM, and PINN-type methods
and then discuss existing bilevel formulations.

In Kansa methods, the kernel family and the lengthscale
\(\ell(\boldsymbol\omega)\) are
central solver-design choices and substantially affect both accuracy and
conditioning \cite{fornberg2015solving,wright2017stable}.
Related developments include leave-one-out cross-validation for lengthscale
selection \cite{rippa1999algorithm}, stable evaluation in the flat-kernel regime
\cite{wright2017stable}, and hybrid Gaussian--cubic RBFs
\cite{mishra2018hybrid}. More recently, neural-network and data-driven methods
have been developed to adapt the lengthscale~\cite{mojarrad2023new,veiga2026learning}.

In RFM/ELM-type solvers, random hidden features
\(\{\psi_j(\cdot;\boldsymbol\omega)\}_{j=1}^{M}\) are held fixed
while the output coefficients \(\boldsymbol\beta\) are fitted by least squares
\cite{dwivedi2020physics,dong2021local}. Approximation accuracy depends
strongly on the parameters that govern feature generation. Existing studies
have used a differential-evolution algorithm \cite{dong2022computing} or
Bayesian optimization \cite{dwivedi2025kernel} to select these parameters.
The corresponding fixed-feature least-squares problems are often severely
ill-conditioned \cite{tan2026high}, making careful selection of the
regularization weight important
\cite{naik2019iterative,yildirim2019performance,naik2020automatic}.

For PINN-type methods, several approaches adjust loss weights during
training, including self-adaptive weighting \cite{mcclenny2023self,xiang2022self},
learning-rate annealing based on gradient statistics
\cite{wang2021understanding}, calibration based on the neural tangent kernel
\cite{wang2022and}, and causality-based weighting
\cite{wang2024respecting}.
These studies motivate including residual weights among the solver-design
parameters \(\boldsymbol\omega\). A separate line of work focuses on learning
the trial representation.
In such cases, \(\boldsymbol\omega\) may include the parameters governing
the trial basis. Examples include
the Gaussian centers and scales optimized in Physics-informed Gaussians
\cite{kang2025pig} and the Fourier-enhanced feature basis updated in IFeF-PINN
\cite{wu2026iterative}.

In the approaches mentioned above, coefficient fitting and
solver-design parameter selection are typically coupled, and the same
collocation points may be used for both tasks, which can increase the risk of
overfitting in parameter selection. A bilevel
formulation makes the roles of these two tasks explicit.
For each choice of \(\boldsymbol\omega\), the lower level fits the
coefficients \(\boldsymbol\beta\),
while the upper level evaluates the resulting approximate solution to select \(\boldsymbol\omega\).

Related bilevel formulations have already been developed for different
mesh-free PDE solvers.
IFeF-PINN \cite{wu2026iterative} and the Gaussian-process bilevel method
of Nelsen et al.~\cite{nelsen2025bilevel} use different strategies for the upper-level update
in nonlinear PDE problems. IFeF-PINN updates the feature parameters using a
partial derivative with the fitted coefficients \(\boldsymbol\beta\) held
fixed. This update omits the hypergradient contribution arising from the
dependence of \(\boldsymbol\beta\) on \(\boldsymbol\omega\).
For example, if a coefficient-regularization weight enters only the lower-level problem, the partial derivative of the upper-level objective with respect to this weight is zero. Consequently, the fixed-coefficient update leaves this weight unchanged, potentially limiting the accuracy of the resulting PDE approximation.
The Gaussian-process bilevel method~\cite{nelsen2025bilevel} instead applies a Gauss--Newton (GN)
linearization to both training and validation residuals and
optimizes the resulting approximate bilevel problem.
In particular, the lower-level GN linearization omits the second-order term in the exact Hessian, which enters the adjoint equation for the exact hypergradient. Omitting this term can introduce large hypergradient errors and thus impair solver-design parameter selection.
Thus, both approaches simplify the hypergradient computation by omitting some key contributions to the hypergradient, but neither quantifies the
resulting hypergradient error or its effect on the final PDE solution error.

Another bilevel approach is the self-optimization
physics-informed Fourier-features randomized neural network (SO-PIFRNN) method
\cite{linghu2025so}, which uses derivative-free particle-swarm optimization to search over
solver-design parameters \(\boldsymbol\omega\). Its upper-level criterion
is the relative \(L^2\) solution error, which requires an analytic or
sufficiently accurate reference solution. Moreover, the search neither uses
upper-level derivative information to guide the parameter updates nor reuses
fitted coefficients \(\boldsymbol\beta\) from earlier evaluations to
initialize subsequent lower-level solves. Consequently, it typically requires many parameter evaluations and repeated lower-level solves, leading to substantial computational cost.

Overall, the aforementioned parameter-selection strategies and
bilevel PDE-solver formulations are limited in three respects.
First, they are tailored to particular solvers and therefore do not provide a
common algorithmic framework applicable across different mesh-free solver
classes. Second, for nonlinear PDEs, existing bilevel methods either rely on costly derivative-free search or simplify the bilevel solve by omitting key contributions to the hypergradient at the expense of solution accuracy. Third, existing studies neither quantify the hypergradient error nor analyze how it affects the PDE solution error.

\subsection{Our Approach and Contributions}
\label{sec:validation-driven-bilevel-formulation}

To address the aforementioned limitations, we introduce a unified validation-driven bilevel optimization framework
for the mesh-free solver classes described in
Section~\ref{sec:residual-coefficient-fitting}. The lower
level fits the coefficients \(\boldsymbol\beta\) for
fixed solver-design parameters \(\boldsymbol\omega\)
using training collocation residuals, while the upper level selects
\(\boldsymbol\omega\) using the residuals of the resulting approximate solution
on held-out validation points.

Let \(\mathcal T=(\mathcal T_I,\mathcal T_B)\) and
\(\mathcal V=(\mathcal V_I,\mathcal V_B)\) denote the separate training and
validation sets of collocation points. For
\(\mathcal C\in\{\mathcal T,\mathcal V\}\),
\(\mathcal C_I\subset\Omega\) and \(\mathcal C_B\subset\partial\Omega\)
are nonempty sets of interior and boundary collocation points, respectively.
We then define the residual vector
\[
    \boldsymbol\rho(\boldsymbol\beta,\boldsymbol\omega;\mathcal C):=
    \begin{bmatrix}
        |\mathcal C_I|^{-1/2}
        \bigl(
            \mathcal L[\widehat{u}(\cdot;\boldsymbol\beta,\boldsymbol\omega)](\boldsymbol x)
            -f(\boldsymbol x)
        \bigr)_{\boldsymbol x\in\mathcal C_I}\\[0.4em]
        |\mathcal C_B|^{-1/2}
        \bigl(
            \mathcal B[\widehat{u}(\cdot;\boldsymbol\beta,\boldsymbol\omega)](\boldsymbol x)
            -g(\boldsymbol x)
        \bigr)_{\boldsymbol x\in\mathcal C_B}
\end{bmatrix}.
\]
With a training-residual weighting matrix \(\mathbf W_{\mathcal T}(\boldsymbol\omega)\)
and the weighted training residual
\(\boldsymbol\rho_{\mathcal T}^{\mathrm{tr}}(\boldsymbol\beta,\boldsymbol\omega)
:=\mathbf W_{\mathcal T}(\boldsymbol\omega)\boldsymbol\rho(\boldsymbol\beta,\boldsymbol\omega;\mathcal T)\),
we define the lower-level objective
\begin{equation*}
    \mathcal E(\boldsymbol\beta,\boldsymbol\omega;\mathcal T):=
    \frac{1}{2}
    \bigl\|
        \boldsymbol\rho_{\mathcal T}^{\rm tr}
        (\boldsymbol\beta,\boldsymbol\omega)
    \bigr\|_2^2
    +
    \frac{\lambda}{2}
    \boldsymbol\beta^\top
    \mathbf R_{\mathcal T}(\boldsymbol\omega)
    \boldsymbol\beta,
\end{equation*}
where \(\lambda\geq 0\) is a regularization parameter and
\(\mathbf R_{\mathcal T}(\boldsymbol\omega)\in\mathbb R^{M\times M}\) is a
symmetric positive-semidefinite regularization matrix. If \(\lambda\) is also
selected at the upper level, we write \(\lambda=\lambda(\boldsymbol\omega)\).

On the validation set, we use the residual
\(\boldsymbol\rho(\boldsymbol\beta,\boldsymbol\omega;\mathcal V)\)
to define the upper-level objective
\begin{equation*}
    \mathcal J(\boldsymbol\beta,\boldsymbol\omega;\mathcal V):=
    \frac{1}{2}
    \bigl\|
        \boldsymbol\rho(\boldsymbol\beta,\boldsymbol\omega;\mathcal V)
    \bigr\|_2^2.
\end{equation*}
The idealized bilevel problem is
\begin{subequations}
\begin{align}
    &\min_{\boldsymbol\omega}\quad
    \mathcal J(
        \boldsymbol\beta^\star(\boldsymbol\omega),
        \boldsymbol\omega;
        \mathcal V)
    \label{eq:outer_objective}\\
    &\text{s.t.}\quad
    \boldsymbol\beta^\star(\boldsymbol\omega)
    \in
    \operatorname*{arg\,min}_{\boldsymbol\beta\in\mathbb R^M}
    \mathcal E(\boldsymbol\beta,\boldsymbol\omega;\mathcal T).
    \label{eq:inner_constraint}
\end{align}
\end{subequations}
In \eqref{eq:inner_constraint}, the lower-level response is required
to be a global minimizer. For nonlinear PDEs, however, the lower-level
objective is generally nonconvex, and we do not assume that the
lower-level solver returns a local or global minimizer. We therefore
analyze a selected branch of lower-level stationary points, also
denoted by \(\boldsymbol\beta^\star(\boldsymbol\omega)\), satisfying
\[
\nabla_{\boldsymbol\beta}\mathcal E(
\boldsymbol\beta^\star(\boldsymbol\omega),
\boldsymbol\omega;\mathcal T)=\mathbf0.
\]
Its local smoothness follows under the Hessian nonsingularity
condition stated in Section~\ref{subsec:hypergradient-bias}.
All subsequent reduced objectives and hypergradients refer to this
selected stationary branch. Substituting this response into the
validation objective gives the reduced objective
\begin{equation}
\label{eq:reduced-validation-objective}
\Phi(\boldsymbol\omega)
:=
\mathcal J(\boldsymbol\beta^\star(\boldsymbol\omega),
\boldsymbol\omega;\mathcal V).
\end{equation}

A gradient-based update of \(\boldsymbol\omega\) requires the hypergradient
\(\nabla\Phi(\boldsymbol\omega)\), which depends on the lower-level coefficient
response \(\boldsymbol\beta^\star(\boldsymbol\omega)\) and the solution of an
associated adjoint equation.
For linear PDEs, fitting the coefficients at fixed
\(\boldsymbol\omega\) reduces to a linear least-squares problem, so both
\(\boldsymbol\beta^\star(\boldsymbol\omega)\) and the adjoint require only a
single linear solve. For nonlinear
PDEs, however, computing the hypergradient presents two main difficulties:
solving the nonlinear lower-level problem to high accuracy is
costly, and the lower-level Hessian in the adjoint system, containing a second-order term
(see~\eqref{eq:npde:local-hessian-decomposition}), is typically
ill-conditioned and indefinite at an approximate response.
In practice, neither the lower-level problem nor the adjoint system needs to be
solved to high accuracy at every upper-level iteration: inexact solves suffice
to form approximate hypergradients for updating
\(\boldsymbol\omega\)~\cite{pedregosa2016hyperparameter,grazzi2020iteration}.

The key task is therefore to exploit the nonlinear least-squares structure to construct
efficient finite-step approximations and quantify how these approximations affect
hypergradient and PDE solution accuracy. Specifically, we approximate the
lower-level response
using finite, warm-started damped GN updates. For the adjoint solve, we use the damped GN matrix to compute an initial
adjoint and to precondition a finite number of minimum residual (MINRES) iterations, fully
leveraging the efficiency of the GN approximation while retaining the
second-order information to improve solution
accuracy. We then analyze the resulting hypergradient error locally and combine this
analysis with validation sampling bounds and residual stability to derive error
bounds for the computed PDE solution. Overall, our main contributions are as follows.
\begin{enumerate}
\item We formulate a unified validation-driven bilevel optimization
framework that decouples solver-design parameter selection from coefficient fitting
in mesh-free PDE solvers, including Kansa, RFM/ELM-type, and
PINN-type methods.

\item We develop a practical and efficient implicit hypergradient algorithm
for nonlinear PDEs that avoids costly high-accuracy lower-level solves
through finite warm-started damped GN updates. For the adjoint solve, we retain the second-order term in the target adjoint system and use the damped
GN matrix to compute an initial adjoint and to precondition MINRES.


\item We derive local hypergradient error bounds for the proposed
algorithm along selected nondegenerate stationary branches,
separating the contributions of the finite lower-level solve,
adjoint damping, and the final adjoint residual.


\item We relate validation residuals to PDE solution accuracy and establish
error bounds for the computed solution in terms of the
best attainable population residual along the selected stationary branch,
validation sampling error, hypergradient error, and lower-level
stationarity error.

\item Numerical experiments across representative mesh-free solvers on several nonlinear PDEs, including a high-dimensional problem, demonstrate the solution accuracy and computational efficiency of the proposed method and illustrate the relation between validation residuals and solution errors.

\end{enumerate}

\subsection{Organization}
The remainder of this paper is organized as follows.
In Section~\ref{sec:algorithmic-design}, we develop the approximate
implicit hypergradient method, including the approximate lower-level
response and the GN-preconditioned MINRES adjoint computation.
The accuracy of the resulting hypergradient approximations is analyzed
in Section~\ref{sec:bias}. Section~\ref{sec:validation-residual-analysis}
then establishes the connection between validation residuals and PDE
solution accuracy and derives error bounds for the computed approximations.
Numerical experiments across representative mesh-free solvers are
presented in Section~\ref{sec:numerical-experiments}.
Finally, Section~\ref{sec:conclusions} summarizes the main findings
and discusses directions for future research.
\section{Algorithm}
\label{sec:algorithmic-design}


In this section, we develop an approximate implicit hypergradient
method for the reduced validation objective
\eqref{eq:reduced-validation-objective} along the selected
lower-level stationary branch.
For notational simplicity, we write
\(\mathcal E(\boldsymbol\beta,\boldsymbol\omega)\equiv\mathcal E(\boldsymbol\beta,\boldsymbol\omega;\mathcal T)\),
\(\mathcal J(\boldsymbol\beta,\boldsymbol\omega)\equiv\mathcal J(\boldsymbol\beta,\boldsymbol\omega;\mathcal V)\), and
\(\boldsymbol\rho(\boldsymbol\beta,\boldsymbol\omega)\equiv\boldsymbol\rho_{\mathcal T}^{\rm tr}(\boldsymbol\beta,\boldsymbol\omega)\in\mathbb R^{N_\rho}\).

In terms of the reduced validation objective \(\Phi(\boldsymbol\omega)\) in
\eqref{eq:reduced-validation-objective}, the basic implicit hypergradient
update~\cite{lorraine2020optimizing,pedregosa2016hyperparameter} reads
\begin{equation}
\label{eq:exact-update}
    \boldsymbol\omega^{k+1}
    =
    \boldsymbol\omega^k-\alpha_k\nabla\Phi(\boldsymbol\omega^k),
    \qquad \alpha_k>0.
\end{equation}%
Under the local smoothness and Hessian nonsingularity assumptions
in Section~\ref{subsec:hypergradient-bias}, differentiating the
stationarity equation and applying the chain rule to
\(\Phi(\boldsymbol\omega)
=\mathcal J(\boldsymbol\beta^\star(\boldsymbol\omega),
\boldsymbol\omega)\)
yields the implicit hypergradient
%
\begin{equation}\label{eq:implicit_gradient}
    \nabla\Phi(\boldsymbol\omega)
    =
    \nabla_{\boldsymbol\omega}
    \mathcal J(
        \boldsymbol\beta^\star(\boldsymbol\omega),
        \boldsymbol\omega)
    -
    \nabla_{\boldsymbol\omega\boldsymbol\beta}^{2}
    \mathcal E(
        \boldsymbol\beta^\star(\boldsymbol\omega),
        \boldsymbol\omega)
    \mathbf v^\star(\boldsymbol\omega),
\end{equation}
where the adjoint \(\mathbf v^\star(\boldsymbol\omega)\) solves the linear system
\begin{equation}
   \label{eq:adjoint}
    \nabla_{\boldsymbol\beta\boldsymbol\beta}^{2}
    \mathcal E(
        \boldsymbol\beta^\star(\boldsymbol\omega),
        \boldsymbol\omega)
    \mathbf v^\star(\boldsymbol\omega)
    =
    \nabla_{\boldsymbol\beta}
    \mathcal J(
        \boldsymbol\beta^\star(\boldsymbol\omega),
        \boldsymbol\omega),
\end{equation}
and the selected lower-level response branch
\(\boldsymbol\beta^\star(\boldsymbol\omega)\) is characterized locally by the
stationarity condition
\begin{equation}
\label{eq:lower-stationarity}
    \nabla_{\boldsymbol\beta}\mathcal E(
        \boldsymbol\beta^\star(\boldsymbol\omega),
        \boldsymbol\omega)=\mathbf0.
\end{equation}
In \eqref{eq:implicit_gradient}, the first term differentiates the validation
objective with the fitted coefficients \(\boldsymbol\beta^\star(\boldsymbol\omega)\)
held fixed, whereas the second term corrects this derivative for the fact that
changing \(\boldsymbol\omega\) also changes the response
\(\boldsymbol\beta^\star(\boldsymbol\omega)\).

Equations~\eqref{eq:exact-update}--\eqref{eq:lower-stationarity}
therefore give a conceptually straightforward exact update, but its implementation for nonlinear PDEs presents two computational bottlenecks.
First, the lower-level stationarity condition
\eqref{eq:lower-stationarity} admits no closed-form solution, and solving it
to high accuracy at every upper-level iteration would dominate the overall cost.
Second,  the adjoint equation \eqref{eq:adjoint} involves
the lower-level Hessian 
\(\nabla_{\boldsymbol\beta\boldsymbol\beta}^{2}\mathcal E\),  which is typically ill-conditioned and indefinite at an approximate lower-level
response, making the adjoint system difficult to solve.
To address both computational bottlenecks, we retain the original nonlinear
objectives and approximate only the two solves in
 \eqref{eq:adjoint} and \eqref{eq:lower-stationarity}.


\subsection{Lower-Level Response Computation via Finite Damped GN}
\label{subsec:algorithm-lower-response}

We first solve the stationarity condition \eqref{eq:lower-stationarity}
inexactly to obtain an approximation to \(\boldsymbol\beta^\star(\boldsymbol\omega)\).
Given the nonlinear least-squares structure of \(\mathcal E\), the GN method is
particularly well suited~\cite{gill1978algorithms}: each step uses only the
residual and its Jacobian, avoiding the expensive computation of the exact
Hessian.

At upper-level iteration \(k\), starting from \(\boldsymbol\beta^{k,0}\), we
take at most \(T\) damped GN steps. At inner iteration
\(t\in\{0,\ldots,T-1\}\), define
\(\boldsymbol\rho^{k,t}:=\boldsymbol\rho(\boldsymbol\beta^{k,t},\boldsymbol\omega^k)\)
and \(\mathbf J^{k,t}:=\nabla_{\boldsymbol\beta}\boldsymbol\rho(\boldsymbol\beta^{k,t},\boldsymbol\omega^k)\).
The lower-level gradient and the damped GN matrix are respectively given by
\begin{equation*}
\begin{aligned}
\boldsymbol g_\beta^{k,t}
&:=
\nabla_{\boldsymbol\beta}
\mathcal E(\boldsymbol\beta^{k,t},\boldsymbol\omega^k)
=
(\mathbf J^{k,t})^\top\boldsymbol\rho^{k,t}
+
\lambda\mathbf R_{\mathcal T}(\boldsymbol\omega^k)
\boldsymbol\beta^{k,t},
\\
\mathbf B^{k,t}
&:=
(\mathbf J^{k,t})^\top\mathbf J^{k,t}
+
\lambda\mathbf R_{\mathcal T}(\boldsymbol\omega^k)
+
\mu_{k,t}\mathbf I.
\end{aligned}
\end{equation*}
Here \(\mu_{k,t}\ge0\) is a damping parameter chosen so that
\(\mathbf B^{k,t}\succeq\underline b_{\mathrm{in}}\mathbf I\) for some
\(\underline b_{\mathrm{in}}>0\), ensuring a well-posed damped GN step:
\begin{equation}
\label{eq:npde:gn-inner-step}
\mathbf B^{k,t}\mathbf p^{k,t}
=-\boldsymbol g_\beta^{k,t},
\qquad
\boldsymbol\beta^{k,t+1}
=
\boldsymbol\beta^{k,t}+\mathbf p^{k,t}.
\end{equation}
Solving every lower-level problem to high accuracy would incur substantial
computational cost. We therefore terminate the lower-level loop either when
\(\|\boldsymbol g_\beta^{k,t}\|_2\le\varepsilon_k\) or when the iteration
budget \(T\) is exhausted, returning the approximate response
\(\widehat{\boldsymbol\beta}^k\). This balances lower-level computational cost and response accuracy,
substantially accelerating the overall bilevel optimization.

To further reduce the lower-level computational cost, we use warm starts
across upper-level iterations. Along the locally smooth response branch,
nearby solver-design parameters have nearby stationary responses.
The response computed at the previous upper-level iteration therefore
provides a natural initialization for the current lower-level
solve~\cite{arbel2022amortized,ji2021bilevel}, and we set
\(\boldsymbol\beta^{k,0}:=\widehat{\boldsymbol\beta}^{k-1}\) for \(k\ge1\).

\subsection{Approximate Hypergradient via GN-Preconditioned MINRES}
\label{subsec:approximate-implicit-gradient}
\label{subsec:algorithm-adjoint-refinement}

We next address the second computational bottleneck: solving the adjoint equation~\eqref{eq:adjoint} at the computed pair \((\widehat{\boldsymbol\beta}^k,\boldsymbol\omega^k)\). Define
\(\mathbf H^k:=\nabla_{\boldsymbol\beta\boldsymbol\beta}^{2}\mathcal E(\widehat{\boldsymbol\beta}^{k},\boldsymbol\omega^k)\)
and \(\boldsymbol b^k:=\nabla_{\boldsymbol\beta}\mathcal J(\widehat{\boldsymbol\beta}^{k},\boldsymbol\omega^k)\),
so that the adjoint system reads \(\mathbf H^k\mathbf v=\boldsymbol b^k\).
With \(\boldsymbol\rho^k:=\boldsymbol\rho(\widehat{\boldsymbol\beta}^{k},\boldsymbol\omega^k)\)
and \(\mathbf J^k:=\nabla_{\boldsymbol\beta}\boldsymbol\rho(\widehat{\boldsymbol\beta}^{k},\boldsymbol\omega^k)\),
the Hessian \(\mathbf H^k\) decomposes as
\begin{equation}
\label{eq:npde:local-hessian-decomposition}
\begin{aligned}
\mathbf H^k
&=(\mathbf J^k)^\top\mathbf J^k
+\lambda\mathbf R_{\mathcal T}(\boldsymbol\omega^k)
+\Delta_{\mathrm{GN}}^k,\\
\Delta_{\mathrm{GN}}^k
&:=\sum_{i=1}^{N_\rho}\rho_i^k
\nabla^2_{\boldsymbol\beta\boldsymbol\beta}
\rho_i(\widehat{\boldsymbol\beta}^{k},\boldsymbol\omega^k).
\end{aligned}
\end{equation}
In~\eqref{eq:npde:local-hessian-decomposition}, the first two terms form the
usual regularized GN matrix, which is symmetric positive semidefinite.
The residual-weighted second-order term \(\Delta_{\mathrm{GN}}^k\) is also
symmetric, but may introduce negative curvature. Consequently, the Hessian
\(\mathbf H^k\) is symmetric but typically indefinite and ill-conditioned.
Standard approximate implicit-differentiation schemes based on the conjugate
gradient (CG) method are designed for symmetric positive-definite Hessian
systems, while fixed-point and Neumann-series approximations require a
contractive iteration~\cite{lorraine2020optimizing,grazzi2020iteration,
ji2021bilevel}; neither condition is guaranteed here.

A tempting way to simplify the adjoint solve is to discard
\(\Delta_{\mathrm{GN}}^k\) and use the GN approximation
\(\mathbf H^k\approx(\mathbf J^k)^\top\mathbf J^k+\lambda\mathbf R_{\mathcal T}(\boldsymbol\omega^k)\).
Because \(\Delta_{\mathrm{GN}}^k\) scales with the residual
\(\boldsymbol\rho^k\), a sufficiently small \(\|\boldsymbol\rho^k\|_2\)
ensures that \(\Delta_{\mathrm{GN}}^k\) is small.
Across the upper-level iterations, however, \(\|\boldsymbol\rho^k\|_2\)
need not be small at \(\widehat{\boldsymbol\beta}^k\), and omitting
\(\Delta_{\mathrm{GN}}^k\) from the adjoint solve can introduce substantial
 error, as analyzed in Remark~\ref{rem:adjoint-refinement-error}.

Rather than omitting \(\Delta_{\mathrm{GN}}^k\) from the adjoint
system, we use a damped GN matrix only to compute an inexpensive initial adjoint
and to precondition MINRES.
Specifically, we define the damped GN matrix
\begin{equation}
\label{eq:npde:gn-hg-matrix}
\mathbf B^k
:=
(\mathbf J^k)^\top\mathbf J^k
+
\lambda\mathbf R_{\mathcal T}(\boldsymbol\omega^k)
+
\mu\mathbf I,
\qquad \mu\ge0,
\end{equation}
with \(\mu\) chosen so that \(\mathbf B^k\) is positive definite, and solve
\(\mathbf B^k\mathbf v_0^k=\boldsymbol b^k\) to obtain \(\mathbf v_0^k\).

To retain \(\Delta_{\mathrm{GN}}^k\), we form the
target adjoint system \(\overline{\mathbf H}^k\mathbf v=\boldsymbol b^k\) using the damped exact-Hessian matrix
\begin{equation}
\label{eq:npde:damped-hessian}
\overline{\mathbf H}^{k}
:=
\mathbf H^k+\mu\mathbf I
=\mathbf B^k+\Delta_{\mathrm{GN}}^k.
\end{equation}
By \(\mathbf B^k\mathbf v_0^k=\boldsymbol b^k\) and
\eqref{eq:npde:damped-hessian}, the initial residual of \(\mathbf v_0^k\) in
this target system is
\(\boldsymbol b^k-\overline{\mathbf H}^{k}\mathbf v_0^k
=-\Delta_{\mathrm{GN}}^k\mathbf v_0^k\), which can be evaluated through a
Hessian--vector product without assembling \(\mathbf H^k\). Given a residual
tolerance \(\tau_{\mathrm{adj}}>0\) and an iteration count \(s\in\mathbb N_0\),
we first check whether
\[
\|\boldsymbol b^k-\overline{\mathbf H}^{k}\mathbf v_0^k\|_2
\le\tau_{\mathrm{adj}}\max\{1,\|\boldsymbol b^k\|_2\}.
\]
If this condition holds, \(\mathbf v_0^k\) already satisfies the prescribed
residual tolerance for the target system, and we set
\(\mathbf v_s^k:=\mathbf v_0^k\) without further iteration. Otherwise, we apply
\(s\) iterations of \(\mathbf B^k\)-preconditioned
MINRES~\cite{paige1975solution} to
\(\overline{\mathbf H}^{k}\mathbf v=\boldsymbol b^k\), starting from
\(\mathbf v_0^k\), to obtain \(\mathbf v_s^k\).

The returned adjoint approximation \(\mathbf v_s^k\) then defines the approximate implicit
hypergradient direction
\begin{equation}
\label{eq:npde:gn-hg-refined}
\boldsymbol d_s^k
:=
\nabla_{\boldsymbol\omega}
\mathcal J(\widehat{\boldsymbol\beta}^{k},\boldsymbol\omega^k)
-
\nabla_{\boldsymbol\omega\boldsymbol\beta}^{2}
\mathcal E(\widehat{\boldsymbol\beta}^{k},\boldsymbol\omega^k)
\mathbf v_s^k.
\end{equation}
Algorithm~\ref{alg:npde-implicit-hypergradient} summarizes this
adjoint and hypergradient computation.

\begin{algorithm}[H]
\caption{Approximate implicit hypergradient via GN-preconditioned MINRES}
\label{alg:npde-implicit-hypergradient}
\begin{algorithmic}[1]
\Require Approximate lower-level response
\(\widehat{\boldsymbol\beta}^{k}\), solver-design parameters \(\boldsymbol\omega^k\), and
adjoint parameters \((\mu,\tau_{\mathrm{adj}},s)\).
\State Solve \(\mathbf B^k\mathbf v_0^k=\boldsymbol b^k\).
\State Set \(\mathbf v_s^k\gets\mathbf v_0^k\).
\If{ \(\|\boldsymbol b^k-\overline{\mathbf H}^{k}\mathbf v_0^k\|_2/
\max\{1,\|\boldsymbol b^k\|_2\}>\tau_{\mathrm{adj}}\)}
    \State Compute \(\mathbf v_s^k\) by applying \(s\) iterations of \(\mathbf B^k\)-preconditioned MINRES to
    \(\overline{\mathbf H}^{k}\mathbf v=\boldsymbol b^k\), starting from \(\mathbf v_0^k\).
\EndIf
\State Compute \(\boldsymbol d_s^k\) by
\eqref{eq:npde:gn-hg-refined}.
\State \Return \((\mathbf v_s^k,\boldsymbol d_s^k)\).
\end{algorithmic}
\end{algorithm}

\begin{remark}
For \(\mu>0\),
\(\mathbf B^k\succeq\mu\mathbf I\succ0\). When \(\mu=0\), full column rank
of \(\mathbf J^k\), or regularization that is positive on
\(\ker(\mathbf J^k)\), is sufficient for positive definiteness.
Under these conditions, we compute the Cholesky factorization
\(\mathbf B^k=\mathbf L^k(\mathbf L^k)^\top\) once. The GN adjoint
\(\mathbf v_0^k\) is then obtained by two triangular solves, and the same
factorization is reused as the MINRES preconditioner.
Let \(C_{\mathrm{HVP}}\) denote the cost of one Hessian--vector product with
\(\overline{\mathbf H}^{k}\), and let \(C_{\mathrm{prec}}\) denote the cost of
one Cholesky-based preconditioner application. Evaluating the initial residual
\(\boldsymbol b^k-\overline{\mathbf H}^{k}\mathbf v_0^k\) requires one such
product, and if MINRES is invoked, its \(s\) iterations require \(s\) additional
products and at most \(s+1\) preconditioner applications. Hence, after forming and factorizing \(\mathbf B^k\), the
additional computational cost is
$
O\!\left((s+1)C_{\mathrm{HVP}}+(s+1)C_{\mathrm{prec}}\right).
$
Exact second-order information enters only through
matrix--vector products.
\end{remark}

\subsection{Complete Approximate Implicit Hypergradient Method}
\label{subsec:complete-aid-algorithm}

Combining the approximate lower-level response from
Section~\ref{subsec:algorithm-lower-response} with the approximate
hypergradient computation in
Section~\ref{subsec:approximate-implicit-gradient} addresses both computational
bottlenecks of the exact update. Replacing \(\nabla\Phi(\boldsymbol\omega^k)\) in
\eqref{eq:exact-update} by the direction \(\boldsymbol d_s^k\)
returned by Algorithm~\ref{alg:npde-implicit-hypergradient}
gives the upper-level update
\begin{equation}
\label{eq:npde:gn-update}
\boldsymbol\omega^{k+1}
=
\boldsymbol\omega^k-\alpha_k\boldsymbol d_s^k.
\end{equation}
Algorithm~\ref{alg:npde-gn} summarizes the complete method.
Exploiting the nonlinear least-squares structure, the method approximates
the lower-level response using finite, warm-started damped GN updates,
avoiding repeated high-accuracy lower-level solves. For the adjoint solve,
it uses the damped GN matrix to compute an initial adjoint and to precondition
a finite number of MINRES iterations, while retaining the second-order term
\(\Delta_{\mathrm{GN}}^k\) in the target adjoint system to improve
hypergradient accuracy. In Section~\ref{sec:bias},
we quantify how the approximations in the lower-level and adjoint solves
propagate into the hypergradient error.

\begin{algorithm}[H]
\caption{Finite-GN implicit hypergradient method with GN-preconditioned MINRES}
\label{alg:npde-gn}
\begin{algorithmic}[1]
\Require Initial pair
\((\boldsymbol\beta^0,\boldsymbol\omega^0)\); upper-level step sizes
\(\{\alpha_k>0\}_{k=0}^{K-1}\); lower-level tolerances
\(\{\varepsilon_k\ge0\}_{k=0}^{K}\); lower-level iteration
budget \(T\in\mathbb N\) and damping rule; admissible
adjoint parameters
\(\mu\ge0,\ \tau_{\mathrm{adj}}>0,\ s\in\mathbb N_0\); upper-level iteration number
\(K\in\mathbb N\).
\For{\(k=0,\ldots,K-1\)}
    \State Set
    \(\boldsymbol\beta^{k,0}\gets\boldsymbol\beta^0\) if \(k=0\), and
    \(\boldsymbol\beta^{k,0}\gets\widehat{\boldsymbol\beta}^{k-1}\)
    otherwise.
    \State Set
    \(\widehat{\boldsymbol\beta}^{k}\gets\boldsymbol\beta^{k,0}\).
    \For{\(t=0,\ldots,T-1\)}
        \If{\(\|\nabla_{\boldsymbol\beta}\mathcal E(
        \boldsymbol\beta^{k,t},\boldsymbol\omega^k)\|_2
        \le \varepsilon_k\)}
            \State \textbf{break}.
        \EndIf
        \State Update \(\boldsymbol\beta^{k,t+1}\) by
        \eqref{eq:npde:gn-inner-step}.
        \State Set
        \(\widehat{\boldsymbol\beta}^{k}\gets\boldsymbol\beta^{k,t+1}\).
    \EndFor
    \State Compute
    \((\mathbf v_s^k,\boldsymbol d_s^k)\) by
    Algorithm~\ref{alg:npde-implicit-hypergradient} at
    \((\widehat{\boldsymbol\beta}^{k},\boldsymbol\omega^k)\).
    \State Set
    \(\boldsymbol\omega^{k+1}\gets
    \boldsymbol\omega^k-\alpha_k\boldsymbol d_s^k\).
\EndFor
\State Apply the same lower-level loop at \(\boldsymbol\omega^K\), initialized with
\(\widehat{\boldsymbol\beta}^{K-1}\), to obtain \(\widehat{\boldsymbol\beta}^{K}\).
\State \Return \((\widehat{\boldsymbol\beta}^{K},\boldsymbol\omega^K)\).
\end{algorithmic}
\end{algorithm}

\section{Hypergradient Analysis for Bilevel Optimization}
\label{sec:bias}

In this section, we analyze the accuracy of the approximate hypergradient direction computed by Algorithm~\ref{alg:npde-implicit-hypergradient}.
Building on local estimates along the selected response branch, we derive
a bound on the hypergradient error in Theorem~\ref{thm:sc-hypergradient-decomposition}
that separates the contributions of the finite lower-level solve, adjoint
damping, and the final adjoint residual.
This bound serves as a key ingredient in Section~\ref{subsec:end-to-end-transfer} for establishing error bounds for the computed PDE solution.

Throughout Sections~\ref{sec:bias} and~\ref{sec:validation-residual-analysis},
we use the shortened notation introduced in Section~\ref{sec:algorithmic-design}.
Unless otherwise stated, \(\|\cdot\|\) denotes the Euclidean norm for finite-dimensional vectors and the induced \(2\)-norm (spectral norm) for matrices.

\subsection{Pointwise Hypergradient Analysis for Nonconvex Lower-Level Problems}
\label{subsec:hypergradient-bias}

This subsection analyzes the reduced objective \(\Phi\) and its hypergradient pointwise along the selected response branch \(\boldsymbol\beta^\star\).
To this end, we first give local conditions under which the implicit hypergradient
representation \eqref{eq:implicit_gradient}--\eqref{eq:adjoint} is valid.

\begin{assumption}[Selected stationary pair and local smoothness]
\label{ass:sc-g1}
Fix a parameter value
\(\overline{\boldsymbol\omega}\in\mathbb R^{n_\omega}\)
and a selected lower-level stationary point
\(\overline{\boldsymbol\beta}\in\mathbb R^M\) satisfying
\begin{equation*}
	\nabla_{\boldsymbol\beta}\mathcal E(
	\overline{\boldsymbol\beta},\overline{\boldsymbol\omega})=\mathbf0.
\end{equation*}
Assume that there is an open neighborhood
\(\mathcal N_0\subset\mathbb R^M\times\mathbb R^{n_\omega}\) of
\((\overline{\boldsymbol\beta},\overline{\boldsymbol\omega})\) on which
\(\mathcal E\in C^2\) and \(\mathcal J\in C^1\).
\end{assumption}

\begin{assumption}[Pointwise lower-level regularity]
\label{ass:sc-g2}
At the selected pair from Assumption~\ref{ass:sc-g1}, the lower-level Hessian is nonsingular:
\begin{equation}
\label{eq:selected-point-hessian}
\det\nabla^2_{\boldsymbol\beta\boldsymbol\beta}\mathcal E(
\overline{\boldsymbol\beta},\overline{\boldsymbol\omega})\ne0.
\end{equation}
\end{assumption}

The next proposition derives the usual implicit root-map derivative
\cite{pedregosa2016hyperparameter,blondel2022implicit}
as a direct
finite-dimensional consequence of the classical implicit-function
theorem; see \cite[Theorem~9.28]{rudin1976principles} and
\cite[Theorem~1B.1]{dontchev2014implicit}.

\begin{proposition}
\label{prop:sc-response-map}
Under Assumptions~\ref{ass:sc-g1}--\ref{ass:sc-g2}, there are neighborhoods
\(U_{\overline{\boldsymbol\omega}}\) of
\(\overline{\boldsymbol\omega}\) and
\(V_{\overline{\boldsymbol\beta}}\) of
\(\overline{\boldsymbol\beta}\), and a \(C^1\) stationary response
\(\boldsymbol\beta^\star:U_{\overline{\boldsymbol\omega}}
\to V_{\overline{\boldsymbol\beta}}\) such that
\begin{equation}
\label{eq:local-selected-branch}
\boldsymbol\beta^\star(\overline{\boldsymbol\omega})
=\overline{\boldsymbol\beta},
\qquad
\nabla_{\boldsymbol\beta}\mathcal E(
\boldsymbol\beta^\star(\boldsymbol\omega),\boldsymbol\omega)=\mathbf0.
\end{equation}
For every \(\boldsymbol\omega\in U_{\overline{\boldsymbol\omega}}\),
\(\boldsymbol\beta^\star(\boldsymbol\omega)\) is the unique
\(\boldsymbol\beta\in V_{\overline{\boldsymbol\beta}}\) satisfying the
stationarity equation in \eqref{eq:local-selected-branch}.
On \(U_{\overline{\boldsymbol\omega}}\), we take this local branch as the
selected response in \eqref{eq:reduced-validation-objective}.
Consequently, the reduced objective \(\Phi\) in
\eqref{eq:reduced-validation-objective} is \(C^1\) near
\(\overline{\boldsymbol\omega}\), and its hypergradient at \(\overline{\boldsymbol\omega}\) is given by
\eqref{eq:implicit_gradient}--\eqref{eq:adjoint}.
\end{proposition}

\begin{proof}
\(\nabla_{\boldsymbol\beta}\mathcal E\) is \(C^1\), and its derivative with
respect to \(\boldsymbol\beta\) at the selected pair is the nonsingular
Hessian in \eqref{eq:selected-point-hessian}. The classical
implicit-function theorem therefore gives the unique local \(C^1\) response
in \eqref{eq:local-selected-branch}. Differentiating its stationarity
equation at \(\overline{\boldsymbol\omega}\) gives
\[
\nabla^2_{\boldsymbol\beta\boldsymbol\beta}\mathcal E(
\overline{\boldsymbol\beta},\overline{\boldsymbol\omega})
D\boldsymbol\beta^\star(\overline{\boldsymbol\omega})
+\nabla^2_{\boldsymbol\beta\boldsymbol\omega}\mathcal E(
\overline{\boldsymbol\beta},\overline{\boldsymbol\omega})=\mathbf0.
\]
Combining this identity with the chain rule for \(\Phi\),
using symmetry of the lower-level Hessian and equality of the mixed partials,
gives \eqref{eq:implicit_gradient} at \(\overline{\boldsymbol\omega}\);
nonsingularity gives the unique adjoint in \eqref{eq:adjoint}.
\end{proof}

\begin{remark}
The following counterexample shows that, if the Hessian nonsingularity in
Assumption~\ref{ass:sc-g2} is dropped, the hypergradient conclusion of
Proposition~\ref{prop:sc-response-map} can fail.

Let \(\beta,\omega\in\mathbb R\), take
\(\lambda=1/4\) and \(R_{\mathcal T}(\omega)=1\), and define
\begin{equation*}
\mathcal E(\beta,\omega)
:=\frac{1}{2}\!\left[
\left(\beta^2-\frac{5}{8}\right)^2+(\beta-\omega)^2\right]
+\frac{1}{8}\beta^2,
\qquad
\mathcal J(\beta,\omega)
:=\frac{1}{2}(\beta^2+\beta+1)^2.
\end{equation*}
Both objectives are \(C^\infty\), and
\(\partial_\beta\mathcal E(\beta,\omega)=2\beta^3-\omega\). Hence the
coercive lower-level objective has the unique global minimizer
\(
\beta^\star(\omega)=\sqrt[3]{\omega/2}
\).
At \((\overline\beta,\overline\omega)=(0,0)\), stationarity holds but
\(
\partial_{\beta\beta}^2\mathcal E(0,0)=0
\), so Assumption~\ref{ass:sc-g2} fails. Writing
\(x=\sqrt[3]{\omega/2}\), the reduced objective is
\[
\Phi(\omega)
=\frac{1}{2}+x+\frac{3}{2}x^2+\frac{\omega}{2}+\frac{1}{2}x^4.
\]
Its leading nonconstant term is proportional to \(\sqrt[3]{\omega}\), so
the difference quotient at zero is unbounded and the hypergradient does not
exist.
\end{remark}

Proposition~\ref{prop:sc-response-map} gives an exact expression for the
implicit hypergradient using the selected response branch in \eqref{eq:local-selected-branch},
whereas Algorithm~\ref{alg:npde-gn} computes an approximate lower-level response.
To analyze the error of this approximate response and of the resulting hypergradient, we impose the following additional local regularity conditions.

\begin{assumption}[Additional regularity for finite-response estimates]
\label{ass:sc-error-regularity}
There is a radius \(r_{\mathcal N}>0\) such that the coefficient neighborhood
\begin{equation*}
\mathcal N_{\overline{\boldsymbol\omega}}
:=
\left\{\boldsymbol\beta\in\mathbb R^M:
\|\boldsymbol\beta-
\overline{\boldsymbol\beta}\|<r_{\mathcal N}\right\}
\end{equation*}
satisfies
\[
\operatorname{cl}(\mathcal N_{\overline{\boldsymbol\omega}})
\times\{\overline{\boldsymbol\omega}\}
\subset\mathcal N_0,
\]
where \(\mathcal N_0\) is the open smoothness neighborhood in
Assumption~\ref{ass:sc-g1}. The maps
\(\nabla_{\boldsymbol\omega}\mathcal J(
\cdot,\overline{\boldsymbol\omega})\),
\(\nabla_{\boldsymbol\beta}\mathcal J(
\cdot,\overline{\boldsymbol\omega})\),
\(\nabla^2_{\boldsymbol\omega\boldsymbol\beta}\mathcal E(
\cdot,\overline{\boldsymbol\omega})\), and
\(\nabla^2_{\boldsymbol\beta\boldsymbol\beta}\mathcal E(
\cdot,\overline{\boldsymbol\omega})\) are Lipschitz on
\(\mathcal N_{\overline{\boldsymbol\omega}}\). Denote by \(L_H\) a
Lipschitz constant for
\(\boldsymbol\beta\mapsto
\nabla^2_{\boldsymbol\beta\boldsymbol\beta}\mathcal E(
\boldsymbol\beta,\overline{\boldsymbol\omega})\). Each residual component
is \(C^2\) with respect to \(\boldsymbol\beta\).
\end{assumption}

The following estimates
bound the response error in terms of the lower-level stationarity residual and control the Hessian inverse used below.

\begin{proposition}
\label{prop:sc-local-response-estimates}
Under Assumptions~\ref{ass:sc-g1}--\ref{ass:sc-g2} and
Assumption~\ref{ass:sc-error-regularity}, define
\[
m:=\frac{1}{2}\sigma_{\min}\!\left(
\nabla^2_{\boldsymbol\beta\boldsymbol\beta}\mathcal E(
\overline{\boldsymbol\beta},\overline{\boldsymbol\omega})
\right)>0
\qquad\text{and}\qquad
\kappa_H:=m^{-1}.
\]
After decreasing \(r_{\mathcal N}\) so that
\(L_Hr_{\mathcal N}\le m\), for every
\(\boldsymbol\beta\in\mathcal N_{\overline{\boldsymbol\omega}}\),
\begin{equation}
\label{eq:sc:response-error}
m\|\boldsymbol\beta-\overline{\boldsymbol\beta}\|
\le
\|\nabla_{\boldsymbol\beta}\mathcal E(
\boldsymbol\beta,\overline{\boldsymbol\omega})\|,
\end{equation}
and
\begin{equation}
\label{eq:sc:hessian-inverse-bound}
\nabla^2_{\boldsymbol\beta\boldsymbol\beta}\mathcal E(
\boldsymbol\beta,\overline{\boldsymbol\omega})
\text{ is nonsingular},
\qquad
\left\|[\nabla^2_{\boldsymbol\beta\boldsymbol\beta}\mathcal E(
\boldsymbol\beta,\overline{\boldsymbol\omega})]^{-1}\right\|
\le\kappa_H.
\end{equation}
\end{proposition}

\begin{proof}
Put \(\boldsymbol h=\boldsymbol\beta-\overline{\boldsymbol\beta}\) and
\[
\mathbf A(\boldsymbol\beta)
:=\int_0^1
\nabla^2_{\boldsymbol\beta\boldsymbol\beta}\mathcal E(
\overline{\boldsymbol\beta}+t\boldsymbol h,
\overline{\boldsymbol\omega})\,dt.
\]
The fundamental theorem of calculus gives
\(
\nabla_{\boldsymbol\beta}\mathcal E(
\boldsymbol\beta,\overline{\boldsymbol\omega})
=\mathbf A(\boldsymbol\beta)\boldsymbol h
\), while Hessian Lipschitz continuity gives
\(
\|\mathbf A(\boldsymbol\beta)-
\nabla^2_{\boldsymbol\beta\boldsymbol\beta}\mathcal E(
\overline{\boldsymbol\beta},\overline{\boldsymbol\omega})\|
\le (L_H/2)\|\boldsymbol h\|
\).
Since
\(\nabla^2_{\boldsymbol\beta\boldsymbol\beta}\mathcal E(
\overline{\boldsymbol\beta},\overline{\boldsymbol\omega})\)
has smallest singular value \(2m\) and
\(L_H\|\boldsymbol h\|\le m\),
the averaged Hessian \(\mathbf A(\boldsymbol\beta)\)
has smallest singular value at least \(m\); this proves
\eqref{eq:sc:response-error}. The same perturbation
estimate without averaging yields
\(
\sigma_{\min}(\nabla^2_{\boldsymbol\beta\boldsymbol\beta}\mathcal E(
\boldsymbol\beta,\overline{\boldsymbol\omega}))
\ge2m-L_H\|\boldsymbol h\|\ge m
\), which proves \eqref{eq:sc:hessian-inverse-bound} with
\(\kappa_H=m^{-1}\).
\end{proof}

The stationarity residual can be evaluated at the computed response,
so \eqref{eq:sc:response-error} bounds the coefficient error without requiring
the exact response \(\overline{\boldsymbol\beta}=\boldsymbol\beta^\star(\overline{\boldsymbol\omega})\).
Theorem~\ref{thm:sc-hypergradient-decomposition} uses
this estimate together with the inverse Hessian bound
\eqref{eq:sc:hessian-inverse-bound} to bound the hypergradient error.

\subsection{Accuracy Analysis of the Approximate Hypergradient Direction}
\label{subsec:finite-response-hypergradient}

We now analyze the error between the approximate hypergradient direction
\(\boldsymbol d_s^k\) returned by Algorithm~\ref{alg:npde-implicit-hypergradient}
and the exact hypergradient \(\nabla\Phi(\boldsymbol\omega^k)\) at a fixed
upper-level iteration \(k\in\{0,\ldots,K-1\}\).
To apply Proposition~\ref{prop:sc-local-response-estimates} at the \(k\)th
iterate, suppose Assumptions~\ref{ass:sc-g1}, \ref{ass:sc-g2}, and~\ref{ass:sc-error-regularity}
hold at \((\overline{\boldsymbol\beta},\overline{\boldsymbol\omega}):=(\boldsymbol\beta^\star(\boldsymbol\omega^k),\boldsymbol\omega^k)\).
Let \(r_{\mathcal N,k}\), \(m_k\), and \(\kappa_{H,k}\) be the local radius
and constants furnished by Proposition~\ref{prop:sc-local-response-estimates},
and define \(\mathcal N_k:=\{\boldsymbol\beta\in\mathbb R^M:\|\boldsymbol\beta-\boldsymbol\beta^\star(\boldsymbol\omega^k)\|<r_{\mathcal N,k}\}\).

Recall from Section~\ref{subsec:approximate-implicit-gradient} that
\(\mathbf H^k\), \(\boldsymbol b^k\), and the damped target matrix
\(\overline{\mathbf H}^k=\mathbf H^k+\mu\mathbf I\) in
\eqref{eq:npde:damped-hessian} are defined at
\((\widehat{\boldsymbol\beta}^k,\boldsymbol\omega^k)\).
The following assumption specifies the local conditions at the computed response needed for the hypergradient error analysis.

\begin{assumption}[Computed-response and target-system regularity]
\label{ass:sc-g3}
The computed response
\(\widehat{\boldsymbol\beta}^k\) lies in \(\mathcal N_k\). At the computed
pair \((\widehat{\boldsymbol\beta}^k,\boldsymbol\omega^k)\), the target
matrix \(\overline{\mathbf H}^k\) defined in
\eqref{eq:npde:damped-hessian} is nonsingular; set
\(\kappa_{\overline H,k}:=\|(\overline{\mathbf H}^k)^{-1}\|<\infty\).
\end{assumption}

With the stationarity residual \(\eta_k:=\|\nabla_{\boldsymbol\beta}\mathcal E(\widehat{\boldsymbol\beta}^k,\boldsymbol\omega^k)\|\),
Assumption~\ref{ass:sc-g3} and Proposition~\ref{prop:sc-local-response-estimates}
give
\begin{equation}
\label{eq:lower-stationarity-certificate}
\|\widehat{\boldsymbol\beta}^k-
\boldsymbol\beta^\star(\boldsymbol\omega^k)\|
\le\frac{\eta_k}{m_k},
\end{equation}
and ensure that \(\mathbf H^k\) is nonsingular with \(\|(\mathbf H^k)^{-1}\|\le\kappa_{H,k}\).

For the initial and returned adjoints \(\mathbf v_0^k\) and \(\mathbf v_s^k\)
from Algorithm~\ref{alg:npde-implicit-hypergradient}, define their residuals
in the damped target system by
\[
\begin{aligned}
\mathbf r_0^k&:=\boldsymbol b^k-\overline{\mathbf H}^k\mathbf v_0^k
=-\Delta_{\mathrm{GN}}^k\mathbf v_0^k,\\
\mathbf r_s^k&:=\boldsymbol b^k-\overline{\mathbf H}^k\mathbf v_s^k.
\end{aligned}
\]
To separate the error due to the approximate lower-level response from that
of the adjoint computation, define the undamped exact-Hessian reference adjoint
and its associated direction at the computed response by
\[
\mathbf v_H^k:=(\mathbf H^k)^{-1}\boldsymbol b^k,\qquad
\boldsymbol d_H^k
:=\nabla_{\boldsymbol\omega}\mathcal J(
\widehat{\boldsymbol\beta}^k,\boldsymbol\omega^k)
-\nabla^2_{\boldsymbol\omega\boldsymbol\beta}\mathcal E(
\widehat{\boldsymbol\beta}^k,\boldsymbol\omega^k)\mathbf v_H^k.
\]
The following theorem gives the main hypergradient error bound of this section,
separating the contributions of the finite lower-level solve, adjoint damping,
and the final adjoint residual.

\begin{theorem}
\label{thm:sc-hypergradient-decomposition}
Suppose Assumptions~\ref{ass:sc-g1}, \ref{ass:sc-g2}, and
\ref{ass:sc-error-regularity} hold at the selected stationary pair
\((\boldsymbol\beta^\star(\boldsymbol\omega^k),\boldsymbol\omega^k)\), and
Assumption~\ref{ass:sc-g3} holds at the associated computed pair.
There is a constant \(C_{\beta,k}>0\) such that
\begin{equation}
\label{eq:sc:eh-bias}
\|\boldsymbol d_H^k-\nabla\Phi(\boldsymbol\omega^k)\|
\le
C_{\beta,k}
\|\widehat{\boldsymbol\beta}^k-
\boldsymbol\beta^\star(\boldsymbol\omega^k)\|
\le
C_{\beta,k}\frac{\eta_k}{m_k}.
\end{equation}
Moreover, with
\(\mathbf A^k:=
\nabla^2_{\boldsymbol\omega\boldsymbol\beta}\mathcal E(
\widehat{\boldsymbol\beta}^k,\boldsymbol\omega^k)\),
the direction \(\boldsymbol d_s^k\) used by
Algorithm~\ref{alg:npde-gn} to update \(\boldsymbol\omega^k\) through
\eqref{eq:npde:gn-update} satisfies the exact decomposition
\begin{equation}
\label{eq:three-way-hypergradient-split}
\begin{aligned}
\boldsymbol d_s^k-\nabla\Phi(\boldsymbol\omega^k)
=
\underbrace{\boldsymbol d_H^k-\nabla\Phi(\boldsymbol\omega^k)}
_{\text{finite lower-level solve}}+
\underbrace{\mu
\mathbf A^k
(\overline{\mathbf H}^k)^{-1}
(\mathbf H^k)^{-1}\boldsymbol b^k}
_{\text{adjoint damping}}
+
\underbrace{\mathbf A^k
(\overline{\mathbf H}^k)^{-1}\mathbf r_s^k}
_{\text{final adjoint residual}}.
\end{aligned}
\end{equation}
Consequently,
\begin{equation}
\label{eq:sc:refined-main-bias}
\begin{aligned}
\|\boldsymbol d_s^k-\nabla\Phi(\boldsymbol\omega^k)\|
\le
C_{\beta,k}\frac{\eta_k}{m_k}
+
\mu\|\mathbf A^k\|\,
\kappa_{\overline H,k}\kappa_{H,k}\|\boldsymbol b^k\|+
\|\mathbf A^k\|\,
\|(\overline{\mathbf H}^k)^{-1}\mathbf r_s^k\|.
\end{aligned}
\end{equation}
If, in addition, \(\eta_k\le\varepsilon_k\) and
\(\|\mathbf r_s^k\|\le\tau_{\rm adj}\max\{1,\|\boldsymbol b^k\|\}\), then there is a constant \(C_{3,k}>0\) such that
\begin{equation}
\label{eq:sc:refined-bias-via-tolerances}
\|\boldsymbol d_s^k-\nabla\Phi(\boldsymbol\omega^k)\|
\le
C_{3,k}\left(
\frac{\varepsilon_k}{m_k}
+\mu+
\tau_{\rm adj}
\right).
\end{equation}
\end{theorem}

\begin{proof}
At fixed \(\boldsymbol\omega^k\), Assumption~\ref{ass:sc-error-regularity},
\eqref{eq:sc:hessian-inverse-bound}, and the inverse perturbation identity
make the exact-Hessian direction locally Lipschitz in
\(\boldsymbol\beta\in\mathcal N_k\). Its values at
\(\boldsymbol\beta^\star(\boldsymbol\omega^k)\) and
\(\widehat{\boldsymbol\beta}^k\) are respectively
\(\nabla\Phi(\boldsymbol\omega^k)\) and \(\boldsymbol d_H^k\); hence
\eqref{eq:sc:eh-bias} follows from
\eqref{eq:lower-stationarity-certificate}. Now set
\(
\overline{\mathbf v}_H^k
:=(\overline{\mathbf H}^k)^{-1}\boldsymbol b^k
\). The resolvent
identity and the definition of the adjoint residual give
\[
\overline{\mathbf v}_H^k-\mathbf v_H^k
=-\mu(\overline{\mathbf H}^k)^{-1}
(\mathbf H^k)^{-1}\boldsymbol b^k,
\qquad
\mathbf v_s^k-\overline{\mathbf v}_H^k
=-(\overline{\mathbf H}^k)^{-1}\mathbf r_s^k.
\]
Substituting these identities into
\(
\boldsymbol d_s^k-\boldsymbol d_H^k
=-\mathbf A^k(\mathbf v_s^k-\mathbf v_H^k)
\)
gives \eqref{eq:three-way-hypergradient-split}; taking norms gives
\eqref{eq:sc:refined-main-bias}. The adjoint-residual condition and
Assumption~\ref{ass:sc-g3} give
\begin{equation}
\label{eq:adjoint-trigger-energy-control}
\|(\overline{\mathbf H}^k)^{-1}\mathbf r_s^k\|
\le
\kappa_{\overline H,k}\tau_{\rm adj}
\max\{1,\|\boldsymbol b^k\|\}.
\end{equation}
Using
\eqref{eq:adjoint-trigger-energy-control} and
\(\eta_k\le\varepsilon_k\) in \eqref{eq:sc:refined-main-bias} gives
\eqref{eq:sc:refined-bias-via-tolerances} with
\(C_{3,k}:=\max\bigl\{C_{\beta,k},\allowbreak\,
\|\mathbf A^k\|\kappa_{\overline H,k}\kappa_{H,k}\|\boldsymbol b^k\|,\allowbreak\,
\|\mathbf A^k\|\kappa_{\overline H,k}\max\{1,\|\boldsymbol b^k\|\}\bigr\}\).
\end{proof}

Theorem~\ref{thm:sc-hypergradient-decomposition} shows that the
accuracy of \(\boldsymbol d_s^k\) can be controlled without solving the
lower-level problem or the adjoint system exactly: once the lower-level and
adjoint residuals meet their tolerances,
\eqref{eq:sc:refined-bias-via-tolerances} bounds the hypergradient error in
terms of \(\varepsilon_k\), \(\tau_{\mathrm{adj}}\), and the damping \(\mu\).
These bounds are then used in Section~\ref{subsec:end-to-end-transfer} to
quantify how the inexact lower-level and adjoint solves affect the error of the
computed PDE solution across the upper-level iterations.

\begin{remark}
\label{rem:adjoint-refinement-error}
Motivated by the GN linearization of Nelsen et al.~\cite[Sec.~2.1]{nelsen2025bilevel},
we examine the effect of omitting \(\Delta_{\mathrm{GN}}^k\) from the adjoint
system at the computed response. We denote the resulting direction by
\(\boldsymbol d_0^k\), obtained by using \(\mathbf v_0^k\) in place of
\(\mathbf v_s^k\) in \eqref{eq:npde:gn-hg-refined}.
The errors of \(\boldsymbol d_0^k\) and
\(\boldsymbol d_s^k\) then differ only in the adjoint-residual term in
\eqref{eq:three-way-hypergradient-split}.

With \(C_{3,k}\) as in the proof of Theorem~\ref{thm:sc-hypergradient-decomposition},
replacing \(\mathbf r_s^k\) by \(\mathbf r_0^k\) in
\eqref{eq:sc:refined-main-bias} gives
\[
\|\boldsymbol d_0^k-\nabla\Phi(\boldsymbol\omega^k)\|
\le C_{3,k}\left(
    \frac{\eta_k}{m_k}+\mu+
    \frac{\|\mathbf r_0^k\|}
         {\max\{1,\|\boldsymbol b^k\|\}}
\right).
\]
Since \(\mathbf r_0^k=-\Delta_{\mathrm{GN}}^k\mathbf v_0^k\) and
\(\Delta_{\mathrm{GN}}^k\) scales with the fitting residual
\(\boldsymbol\rho^k\), which need not vanish even at an exact lower-level
response, the last term need not vanish as \(\eta_k\to0\) and \(\mu\to0\), and
no tolerance controls it. By contrast,
Algorithm~\ref{alg:npde-implicit-hypergradient} retains
\(\Delta_{\mathrm{GN}}^k\) and uses MINRES to reduce the initial residual
\(\mathbf r_0^k\); when \(\mathbf r_s^k\) meets the residual tolerance, this
term is replaced by \(\tau_{\mathrm{adj}}\) in
\eqref{eq:sc:refined-bias-via-tolerances}.
\end{remark}

\section{Validation Residuals and End-to-End PDE Error Control}
\label{sec:validation-residual-analysis}

Section~\ref{sec:bias} bounds the error of the approximate hypergradient of the
reduced validation objective \(\Phi\); it remains to show how optimizing \(\Phi\) with
inexact solves controls the error of the computed PDE solution. We establish
this control in two steps. First,
Section~\ref{subsec:validation-residual-error-control} relates the validation
objective to the PDE solution error along the selected response branch, and
Section~\ref{subsec:exampleA-verification} verifies the required conditions for
the nonlinear elliptic Kansa solver implemented in Section \ref{sec:numerical-experiments}. Second,
Section~\ref{subsec:end-to-end-transfer} accounts for the finite lower-level
solve and, for the basic
upper-level update, incorporates the hypergradient errors analyzed in
Section~\ref{sec:bias} to obtain the finite-iteration and asymptotic PDE error
bounds.

\subsection{Validation-Based Bilevel Selection and PDE Error Control}
\label{subsec:validation-residual-error-control}

The upper-level objective measures residuals only at the validation points,
whereas PDE stability controls the solution error by a problem-dependent
residual norm. To relate the validation objective to this norm, we introduce
a population residual norm and combine a uniform sampling bound with norm
compatibility, yielding the PDE error bound in Corollary~\ref{cor:local-class-a-posteriori}.
On the selected response branch, Theorem~\ref{thm:bilevel-validation-oracle}
gives a PDE error bound for validation-based selection.

For a trial function \(v\) whose residuals are well defined, set
\(
\mathcal R(v)
:=
\bigl(\mathcal L[v]-f,\,\mathcal B[v]-g\bigr)
\)
and define the problem-dependent residual norm
\begin{equation}
\label{eq:continuous-residual-norm}
\|\mathcal R(v)\|_Y^2
:=
\|\mathcal L[v]-f\|_{Y_I}^2
+
\eta_B\|\mathcal B[v]-g\|_{Y_B}^2,
\qquad \eta_B>0.
\end{equation}
Let \(\nu_I\) and \(\nu_B\) be the probability measures used for interior and
boundary validation sampling, with \(\operatorname{supp}(\nu_I)=\overline\Omega\)
and \(\operatorname{supp}(\nu_B)=\partial\Omega\).
For \(r=(r_I,r_B)\), define the corresponding population residual norm
\begin{equation*}
\|r\|_{Y_0}^2
:=
\int_\Omega |r_I|^2\,d\nu_I
+
\int_{\partial\Omega}|r_B|^2\,d\nu_B.
\end{equation*}
If
\(\mathcal V_I=\{\boldsymbol X_i^I\}_{i=1}^{N_I}\) and
\(\mathcal V_B=\{\boldsymbol X_j^B\}_{j=1}^{N_B}\), with
\(N_I,N_B\ge1\), define the empirical residual norm
\begin{equation}
\label{eq:validation-sampling-operator}
\|\mathsf S_{\mathcal V}r\|^2
:=
\frac{1}{N_I}\sum_{i=1}^{N_I}|r_I(\boldsymbol X_i^I)|^2
+
\frac{1}{N_B}\sum_{j=1}^{N_B}|r_B(\boldsymbol X_j^B)|^2.
\end{equation}
The normalizations in \eqref{eq:validation-sampling-operator} agree with the
validation residual vector \(\boldsymbol\rho(\boldsymbol\beta,\boldsymbol\omega;\mathcal V)\), so
\begin{equation}
\label{eq:validation-loss-sampling-norm}
\|\mathsf S_{\mathcal V}
\mathcal R(\widehat{u}(\cdot;\boldsymbol\beta,\boldsymbol\omega))\|^2
=
2\mathcal J(\boldsymbol\beta,\boldsymbol\omega;\mathcal V).
\end{equation}

Let
\(\mathcal K_{\mathrm{loc}}\subset
\mathbb R^M\times\mathbb R^{n_\omega}\) be
nonempty and compact and specified independently
of the validation sample, and define
\[
\mathcal F_{\mathrm{loc}}
:=
\left\{
\mathcal R(\widehat{u}(\cdot;\boldsymbol\beta,\boldsymbol\omega)):
(\boldsymbol\beta,\boldsymbol\omega)\in\mathcal K_{\mathrm{loc}}
\right\}.
\]
Let \(X\) denote the solution space
in which the PDE error is measured. Population control in \(Y_0\) yields
solution error control in \(X\) if the \(Y_0\)-norm controls the \(Y\)-norm
and residual stability transfers \(Y\)-control to \(X\). The following
assumption states these two links on the fixed candidate class.

\begin{assumption}[Residual stability and norm compatibility]
\label{ass:residual-stability-compatibility}
There are constants \(C_{\mathrm{stab}},C_{\mathrm{loc}}>0\), independent of
\((\boldsymbol\beta,\boldsymbol\omega)\in\mathcal K_{\mathrm{loc}}\), such
that
\begin{align}
\|\widehat{u}(\cdot;\boldsymbol\beta,\boldsymbol\omega)-u\|_X
&\le
C_{\mathrm{stab}}
\|\mathcal R(\widehat{u}(\cdot;\boldsymbol\beta,\boldsymbol\omega))\|_Y,
\label{eq:continuous-residual-stability}
\\
\|r\|_Y
&\le
C_{\mathrm{loc}}\|r\|_{Y_0},
\qquad r\in\mathcal F_{\mathrm{loc}}.
\label{eq:residual-norm-compatibility}
\end{align}
\end{assumption}

Related estimates linking residuals to solution errors appear in
a posteriori analyses for neural PDE solvers
\cite{berrone2022solving,hillebrecht2023rigorous,mukherjee2026rigorous}
and in stability analyses for mesh-free discretizations
\cite{schaback2016all,bohmer2020nonlinear}.

To connect the validation objective in \eqref{eq:validation-loss-sampling-norm}
with the population residual norm, we bound the difference between the empirical
and population squared residual norms uniformly over the residual class \(\mathcal F_{\mathrm{loc}}\).
We use covering numbers to quantify its size. For a subset \(\mathcal A\)
of a metric space with metric \(d\), let \(\mathcal N(\epsilon,\mathcal A;d)\)
denote the minimum cardinality of an internal \(\epsilon\)-net of \(\mathcal A\).
On \(C(\overline\Omega)\times C(\partial\Omega)\), use the uniform metric
\[
d_\infty(r,\widetilde r)
:=
\max\left\{
\|r_I-\widetilde r_I\|_{L^\infty(\Omega)},
\|r_B-\widetilde r_B\|_{L^\infty(\partial\Omega)}
\right\}.
\]
For residual classes, write
\(\mathcal N_\infty(\epsilon,\mathcal F)
:=\mathcal N(\epsilon,\mathcal F;d_\infty)\).

\begin{proposition}
\label{prop:local-class-validation-sampling}
Assume that the residual map \((\boldsymbol\beta,\boldsymbol\omega)\mapsto\mathcal R\!\left(\widehat u(\cdot;\boldsymbol\beta,\boldsymbol\omega)\right)\)
is continuous from
\(\mathcal K_{\mathrm{loc}}\) into
\(C(\overline\Omega)\times C(\partial\Omega)\), equipped with \(d_\infty\).
Set
\[
B_I:=\sup_{r\in\mathcal F_{\mathrm{loc}}}
\|r_I\|_{L^\infty(\Omega)},
\qquad
B_B:=\sup_{r\in\mathcal F_{\mathrm{loc}}}
\|r_B\|_{L^\infty(\partial\Omega)}.
\]
Let the interior and boundary validation points be independent samples from
\(\nu_I\) and \(\nu_B\). For every
\(\epsilon_{\mathrm{cov}}>0\) and \(\zeta\in(0,1)\), define
\begin{equation*}
\delta_{\mathrm{samp}}^2
:=
4\epsilon_{\mathrm{cov}}(B_I+B_B)
+
B_I^2
\sqrt{
\frac{
\log\!\left(
4\mathcal N_\infty(\epsilon_{\mathrm{cov}},
\mathcal F_{\mathrm{loc}})/\zeta
\right)}{2N_I}}+
B_B^2
\sqrt{
\frac{
\log\!\left(
4\mathcal N_\infty(\epsilon_{\mathrm{cov}},
\mathcal F_{\mathrm{loc}})/\zeta
\right)}{2N_B}}.
\end{equation*}
Then, with probability at least \(1-\zeta\),
\begin{equation}
\label{eq:local-class-validation-sampling}
\sup_{r\in\mathcal F_{\mathrm{loc}}}
\left|
\|r\|_{Y_0}^2-
\|\mathsf S_{\mathcal V}r\|^2
\right|
\le
\delta_{\mathrm{samp}}^2.
\end{equation}
\end{proposition}

\begin{proof}
The residual class \(\mathcal F_{\mathrm{loc}}\) is compact in
\(d_\infty\), so a finite internal \(\epsilon_{\mathrm{cov}}\)-net exists and \(B_I,B_B<\infty\). Fix an
internal net
\(\{r^{(m)}\}_{m=1}^{M_{\mathrm{cov}}}\), where
\(M_{\mathrm{cov}}=
\mathcal N_\infty(\epsilon_{\mathrm{cov}},\mathcal F_{\mathrm{loc}})\).
Applying Hoeffding's inequality at each net element and taking a
union bound over the net and both residual components gives,
simultaneously for all \(m\),
\[
\left|
\|r_I^{(m)}\|_{L^2(\nu_I)}^2
-
\frac{1}{N_I}\sum_{i=1}^{N_I}
|r_I^{(m)}(\boldsymbol X_i^I)|^2
\right|
\le t_I,~\text{and}~\left|
\|r_B^{(m)}\|_{L^2(\nu_B)}^2
-
\frac{1}{N_B}\sum_{j=1}^{N_B}
|r_B^{(m)}(\boldsymbol X_j^B)|^2
\right|
\le t_B,
\]
with total failure probability at most \(\zeta\), where
\[
t_I:=B_I^2
\sqrt{\frac{\log(4M_{\mathrm{cov}}/\zeta)}{2N_I}},
\qquad
t_B:=B_B^2
\sqrt{\frac{\log(4M_{\mathrm{cov}}/\zeta)}{2N_B}}.
\]
To extend these bounds to the full residual class, choose for each
\(r\in\mathcal F_{\mathrm{loc}}\) a net element \(r^{(m)}\) with
\(d_\infty(r,r^{(m)})\le\epsilon_{\mathrm{cov}}\).
Since \(r,r^{(m)}\in\mathcal F_{\mathrm{loc}}\),
\[
\bigl||r_I|^2-|r_I^{(m)}|^2\bigr|
\le2B_I\epsilon_{\mathrm{cov}},
\qquad
\bigl||r_B|^2-|r_B^{(m)}|^2\bigr|
\le2B_B\epsilon_{\mathrm{cov}}.
\]
Replacing \(r^{(m)}\) by \(r\) changes each of the population and empirical
squared norms by at most \(2\epsilon_{\mathrm{cov}}(B_I+B_B)\). Hence
\[
\left|
\|r\|_{Y_0}^2-
\|\mathsf S_{\mathcal V}r\|^2
\right|
\le
4\epsilon_{\mathrm{cov}}(B_I+B_B)+t_I+t_B
=
\delta_{\mathrm{samp}}^2.
\]
Taking the supremum proves \eqref{eq:local-class-validation-sampling}.
\end{proof}

Combining this uniform sampling bound with Assumption~\ref{ass:residual-stability-compatibility}
yields the following PDE error bound on \(\mathcal K_{\mathrm{loc}}\).

\begin{corollary}
\label{cor:local-class-a-posteriori}
Under Assumption~\ref{ass:residual-stability-compatibility}, on the event
\eqref{eq:local-class-validation-sampling}, every
\((\boldsymbol\beta,\boldsymbol\omega)\in\mathcal K_{\mathrm{loc}}\) satisfies
\begin{equation*}
\|\widehat{u}(\cdot;\boldsymbol\beta,\boldsymbol\omega)-u\|_X
\le
C_{\mathrm{stab}}C_{\mathrm{loc}}
\left(
2\mathcal J(\boldsymbol\beta,\boldsymbol\omega;\mathcal V)
+
\delta_{\mathrm{samp}}^2
\right)^{1/2}.
\end{equation*}
\end{corollary}

\begin{proof}
For
\(r=\mathcal R(\widehat{u}(\cdot;\boldsymbol\beta,\boldsymbol\omega))\), use
\eqref{eq:continuous-residual-stability},
\eqref{eq:residual-norm-compatibility}, and
\eqref{eq:local-class-validation-sampling} to obtain
\[
\|\widehat{u}-u\|_X
\le
C_{\mathrm{stab}}C_{\mathrm{loc}}
\|r\|_{Y_0}
\le
C_{\mathrm{stab}}C_{\mathrm{loc}}
\left(
\|\mathsf S_{\mathcal V}r\|^2
+
\delta_{\mathrm{samp}}^2
\right)^{1/2}.
\]
Equation~\eqref{eq:validation-loss-sampling-norm} proves the claim.
\end{proof}

We now restrict the candidates to the selected response branch \(\boldsymbol\beta^\star\).
Let \(\mathcal W\subset\mathbb R^{n_\omega}\) be a nonempty compact parameter
set on which \(\boldsymbol\beta^\star\) is defined, with \(\mathcal W\) and
\(\boldsymbol\beta^\star\) fixed independently of the validation sample.
Write \(r_\star(\boldsymbol\omega):=\mathcal R\!\left(\widehat u(\cdot;\boldsymbol\beta^\star(\boldsymbol\omega),\boldsymbol\omega)\right)\)
for the response-residual map, and define the selected lower-level response
graph and the response residual class by
\begin{equation*}
\mathcal G_\star
:=
\left\{
(\boldsymbol\beta^\star(\boldsymbol\omega),\boldsymbol\omega):
\boldsymbol\omega\in\mathcal W
\right\},
\qquad
\mathcal F_\star
:=
\left\{
r_\star(\boldsymbol\omega):
\boldsymbol\omega\in\mathcal W
\right\}.
\end{equation*}

For selection along \(\mathcal G_\star\), uniform sampling control is needed
only over \(\mathcal F_\star\). Applying the argument of Proposition~\ref{prop:local-class-validation-sampling}
to this class yields the following corollary.

\begin{corollary}
\label{cor:bilevel-response-validation-sampling}
Assume that \(r_\star\) is continuous from \(\mathcal W\) into \(C(\overline\Omega)\times C(\partial\Omega)\), equipped with \(d_\infty\). Set
\[
B_{I,\star}
:=
\sup_{r\in\mathcal F_\star}
\|r_I\|_{L^\infty(\Omega)},
\qquad
B_{B,\star}
:=
\sup_{r\in\mathcal F_\star}
\|r_B\|_{L^\infty(\partial\Omega)}.
\]
Let the interior and boundary validation points be independent samples from
\(\nu_I\) and \(\nu_B\). For every
\(\epsilon_{\mathrm{cov}}>0\) and \(\zeta\in(0,1)\), define
\begin{equation}
\label{eq:bilevel-response-class-sampling-floor}
\begin{aligned}
\delta_{\mathrm{samp},\star}^2
:=
4\epsilon_{\mathrm{cov}}(B_{I,\star}+B_{B,\star})
+
B_{I,\star}^2
\sqrt{
\frac{
\log\!\left(
4\mathcal N_\infty(\epsilon_{\mathrm{cov}},
\mathcal F_\star)/\zeta
\right)}{2N_I}}
+
B_{B,\star}^2
\sqrt{
\frac{
\log\!\left(
4\mathcal N_\infty(\epsilon_{\mathrm{cov}},
\mathcal F_\star)/\zeta
\right)}{2N_B}}.
\end{aligned}
\end{equation}
Then, with probability at least \(1-\zeta\),
\begin{equation}
\label{eq:bilevel-response-class-validation-sampling}
\sup_{r\in\mathcal F_\star}
\left|
\|r\|_{Y_0}^2-
\|\mathsf S_{\mathcal V}r\|^2
\right|
\le
\delta_{\mathrm{samp},\star}^2.
\end{equation}
\end{corollary}

\begin{proof}
Since \(\mathcal W\) is compact and \(r_\star\) is continuous,
\(\mathcal F_\star=r_\star(\mathcal W)\) is compact in \(d_\infty\), and
the proof of Proposition~\ref{prop:local-class-validation-sampling} applied
to \(\mathcal F_\star\) yields
\eqref{eq:bilevel-response-class-validation-sampling}.
\end{proof}

Recall the reduced validation objective \(\Phi(\boldsymbol\omega)=\mathcal J(\boldsymbol\beta^\star(\boldsymbol\omega),\boldsymbol\omega;\mathcal V)\)
from \eqref{eq:reduced-validation-objective}, and define its infimum and
the best squared population residual norm along the selected response branch
over \(\mathcal W\) by
\[
\Phi_\star:=\inf_{\boldsymbol\omega\in\mathcal W}\Phi(\boldsymbol\omega),
\qquad
R_\star:=\inf_{\boldsymbol\omega\in\mathcal W}
\left\|\mathcal R\!\left(\widehat{u}(\cdot;\boldsymbol\beta^\star(\boldsymbol\omega),\boldsymbol\omega)\right)\right\|_{Y_0}^2.
\]
Combining Corollary~\ref{cor:bilevel-response-validation-sampling} with Assumption~\ref{ass:residual-stability-compatibility}
then bounds the population residual and PDE solution error at an approximate
minimizer of \(\Phi\) in terms of \(R_\star\).

\begin{theorem}
\label{thm:bilevel-validation-oracle}
Assume the hypotheses of
Corollary~\ref{cor:bilevel-response-validation-sampling} and
Assumption~\ref{ass:residual-stability-compatibility}, and suppose
\(\mathcal G_\star\subset\mathcal K_{\mathrm{loc}}\).
Let
\(\widetilde{\boldsymbol\omega}\in\mathcal W\), possibly selected using the
validation set, satisfy
\begin{equation}
\label{eq:bilevel-approx-upper-optimality}
\Phi(\widetilde{\boldsymbol\omega})
\le
\Phi_\star+\varepsilon_{\mathrm{out}},
\qquad
\varepsilon_{\mathrm{out}}\ge0.
\end{equation}
Then, on the event
\eqref{eq:bilevel-response-class-validation-sampling},
\begin{equation}
\label{eq:bilevel-upper-oracle-comparison}
2\Phi(\widetilde{\boldsymbol\omega})
\le
R_\star+\delta_{\mathrm{samp},\star}^2
+2\varepsilon_{\mathrm{out}},
\end{equation}
and
\begin{equation}
\label{eq:bilevel-continuous-residual-oracle}
\left\|
\mathcal R\!\left(
\widehat{u}(\cdot;
\boldsymbol\beta^\star(\widetilde{\boldsymbol\omega}),
\widetilde{\boldsymbol\omega})
\right)
\right\|_{Y_0}^2
\le
R_\star+2\delta_{\mathrm{samp},\star}^2
+2\varepsilon_{\mathrm{out}}.
\end{equation}
Consequently,
\begin{equation}
\label{eq:bilevel-pde-oracle}
\left\|
\widehat{u}(\cdot;
\boldsymbol\beta^\star(\widetilde{\boldsymbol\omega}),
\widetilde{\boldsymbol\omega})
-u
\right\|_X
\le
C_{\mathrm{stab}}C_{\mathrm{loc}}
\left(
R_\star
+2\delta_{\mathrm{samp},\star}^2
+2\varepsilon_{\mathrm{out}}
\right)^{1/2}.
\end{equation}
If \(\widetilde{\boldsymbol\omega}\) is an exact upper-level minimizer, then one may take
\(\varepsilon_{\mathrm{out}}=0\).
\end{theorem}

\begin{proof}
For any \(\epsilon>0\), choose
\(\boldsymbol\omega_\epsilon\in\mathcal W\) such that
\[
r_\epsilon
:=
\mathcal R\!\left(
\widehat{u}(\cdot;
\boldsymbol\beta^\star(\boldsymbol\omega_\epsilon),
\boldsymbol\omega_\epsilon)
\right)
\]
satisfies
\[
\|r_\epsilon\|_{Y_0}^2
\le
R_\star+\epsilon.
\]
Since
\(\Phi_\star\le\Phi(\boldsymbol\omega_\epsilon)\),
\eqref{eq:bilevel-approx-upper-optimality} and
\eqref{eq:bilevel-response-class-validation-sampling} give
\[
\begin{aligned}
2\Phi(\widetilde{\boldsymbol\omega})
\le
2\Phi(\boldsymbol\omega_\epsilon)
+2\varepsilon_{\mathrm{out}}
=
\|\mathsf S_{\mathcal V}r_\epsilon\|^2
+2\varepsilon_{\mathrm{out}}\le
R_\star+\epsilon+\delta_{\mathrm{samp},\star}^2
+2\varepsilon_{\mathrm{out}}.
\end{aligned}
\]
Letting \(\epsilon\downarrow0\) proves
\eqref{eq:bilevel-upper-oracle-comparison}.
Applying \eqref{eq:bilevel-response-class-validation-sampling} to
\(r_\star(\widetilde{\boldsymbol\omega})\) yields
\[
\left\|
\mathcal R\!\left(
\widehat{u}(\cdot;
\boldsymbol\beta^\star(\widetilde{\boldsymbol\omega}),
\widetilde{\boldsymbol\omega})
\right)
\right\|_{Y_0}^2
\le
2\Phi(\widetilde{\boldsymbol\omega})
+\delta_{\mathrm{samp},\star}^2,
\]
which, together with
\eqref{eq:bilevel-upper-oracle-comparison}, proves
\eqref{eq:bilevel-continuous-residual-oracle}.
Finally, Assumption~\ref{ass:residual-stability-compatibility} gives
\eqref{eq:bilevel-pde-oracle}.
\end{proof}

Theorem~\ref{thm:bilevel-validation-oracle} shows that the reduced
validation objective \(\Phi\) is a computable selection criterion that controls
the PDE solution error along the selected response branch.
Specifically, \eqref{eq:bilevel-pde-oracle} bounds the PDE error in
terms of the best attainable population residual \(R_\star\), the validation
sampling error \(\delta_{\mathrm{samp},\star}^2\), and the upper-level
objective gap \(\varepsilon_{\mathrm{out}}\). Decreasing \(\Phi\) toward
\(\Phi_\star\) reduces \(\varepsilon_{\mathrm{out}}\) and drives this bound
toward the floor
\(C_{\mathrm{stab}}C_{\mathrm{loc}}(R_\star+2\delta_{\mathrm{samp},\star}^2)^{1/2}\).

In Proposition~\ref{prop:local-class-validation-sampling} and
Corollary~\ref{cor:bilevel-response-validation-sampling}, the sampling error
bounds depend on the covering number of the residual class over which uniform
sampling control is required. In the bilevel formulation, the lower level determines
\(\boldsymbol\beta^\star(\boldsymbol\omega)\) using the training set
\(\mathcal T\), so validation-based selection on \(\mathcal V\) proceeds along
this response branch and requires uniform sampling control only over
\(\mathcal F_\star\). For comparison, jointly minimizing
\(\mathcal J(\boldsymbol\beta,\boldsymbol\omega;\mathcal V)\) uses the same
validation points both to fit \(\boldsymbol\beta\) and to select
\(\boldsymbol\omega\). Over a nonempty candidate set
\(\mathcal K_{\mathrm{joint}}\subset\mathbb R^{M+n_\omega}\) fixed independently
of \(\mathcal V\), this requires uniform sampling control over the joint
residual class
\(\mathcal F_{\mathrm{joint}}
:=\left\{
\mathcal R\!\left(\widehat u(\cdot;\boldsymbol\beta,\boldsymbol\omega)\right):
(\boldsymbol\beta,\boldsymbol\omega)\in\mathcal K_{\mathrm{joint}}
\right\}\).
Assume that \(\mathcal W\subset B(0,R_\omega)\subset\mathbb R^{n_\omega}\)
and \(\mathcal K_{\mathrm{joint}}\subset B(0,R_{\mathrm{joint}})\subset\mathbb R^{M+n_\omega}\),
and write \(\mathcal N_2(\epsilon,\mathcal A):=\mathcal N(\epsilon,\mathcal A;\|\cdot\|)\)
for Euclidean subsets \(\mathcal A\). Measuring Lipschitz continuity from the
Euclidean norm to \(d_\infty\), we obtain the following covering bounds for
\(\mathcal F_\star\) and \(\mathcal F_{\mathrm{joint}}\).

\begin{corollary}
\label{cor:bilevel-response-covering}
Assume that \(r_\star\) is \(L_\star\)-Lipschitz on \(\mathcal W\) and the
residual map \((\boldsymbol\beta,\boldsymbol\omega)\mapsto\mathcal R(\widehat u(\cdot;\boldsymbol\beta,\boldsymbol\omega))\)
is \(L_R\)-Lipschitz on \(\mathcal K_{\mathrm{joint}}\), with
\(L_\star,L_R>0\). Then, for every \(\epsilon>0\),
\begin{align}
\mathcal N_\infty(\epsilon,\mathcal F_\star)
&\le\mathcal N_2\!\left(\frac{\epsilon}{L_\star},\mathcal W\right)
\le\left(1+\frac{2R_\omega L_\star}{\epsilon}\right)^{n_\omega},
\label{eq:bilevel-response-covering}\\
\mathcal N_\infty(\epsilon,\mathcal F_{\mathrm{joint}})
&\le\mathcal N_2\!\left(\frac{\epsilon}{L_R},\mathcal K_{\mathrm{joint}}\right)
\le\left(1+\frac{2R_{\mathrm{joint}}L_R}{\epsilon}\right)^{M+n_\omega}.
\label{eq:joint-residual-covering}
\end{align}
If, in addition, the same
\(L_R\)-Lipschitz bound holds on \(\mathcal G_\star\) and
\(\boldsymbol\beta^\star\) is \(L_\beta\)-Lipschitz on \(\mathcal W\), then one may take \(L_\star\le L_R\sqrt{1+L_\beta^2}\).
\end{corollary}

\begin{proof}
If \(\{\boldsymbol\omega_j\}_{j=1}^N\) is an internal
\(\epsilon/L_\star\)-net of \(\mathcal W\) with
\(N=\mathcal N_2(\epsilon/L_\star,\mathcal W)\), then the
\(L_\star\)-Lipschitz property of \(r_\star\) implies
\(\min_j d_\infty(r_\star(\boldsymbol\omega),r_\star(\boldsymbol\omega_j))
\le L_\star\min_j\|\boldsymbol\omega-\boldsymbol\omega_j\|\le\epsilon\)
for every \(\boldsymbol\omega\in\mathcal W\). Thus
\(\{r_\star(\boldsymbol\omega_j)\}_{j=1}^N\) is an internal
\(\epsilon\)-net of \(\mathcal F_\star\), proving the first inequality in
\eqref{eq:bilevel-response-covering}; the same argument on
\(\mathcal K_{\mathrm{joint}}\) gives the first inequality in
\eqref{eq:joint-residual-covering}. The second inequalities in
\eqref{eq:bilevel-response-covering}--\eqref{eq:joint-residual-covering}
follow from the standard Euclidean covering bound
\(\mathcal N_2(\delta,\mathcal A)\le(1+2R/\delta)^d\) for
\(\mathcal A\subset B(0,R)\subset\mathbb R^d\). Finally, under the additional
assumptions, any \(\boldsymbol\omega_1,\boldsymbol\omega_2\in\mathcal W\) satisfy
\[
\begin{aligned}
d_\infty\!\left(
r_\star(\boldsymbol\omega_1),
r_\star(\boldsymbol\omega_2)
\right)
\le
L_R
\left(
\|\boldsymbol\beta^\star(\boldsymbol\omega_1)
-\boldsymbol\beta^\star(\boldsymbol\omega_2)\|^2
+
\|\boldsymbol\omega_1-\boldsymbol\omega_2\|^2
\right)^{1/2}\le
L_R\sqrt{1+L_\beta^2}\,
\|\boldsymbol\omega_1-\boldsymbol\omega_2\|,
\end{aligned}
\]
yielding the stated bound on \(L_\star\).
\end{proof}

Substituting \eqref{eq:bilevel-response-covering} into
\eqref{eq:bilevel-response-class-sampling-floor} gives
\[
\delta_{\mathrm{samp},\star}^2
\le 4\epsilon_{\mathrm{cov}}(B_{I,\star}+B_{B,\star})
+\left(\frac{B_{I,\star}^2}{\sqrt{2N_I}}+
\frac{B_{B,\star}^2}{\sqrt{2N_B}}\right)
\left(\log\frac{4}{\zeta}
+n_\omega\log\!\left(1+\frac{2R_\omega L_\star}{\epsilon_{\mathrm{cov}}}\right)\right)^{1/2},
\]
which makes the dependence of the bounds in
Theorem~\ref{thm:bilevel-validation-oracle} on \(n_\omega\) explicit.
By contrast, applying the argument of
Proposition~\ref{prop:local-class-validation-sampling} to
\(\mathcal F_{\mathrm{joint}}\) together with \eqref{eq:joint-residual-covering}
yields a bound of the same form with dimension factor \(M+n_\omega\) in place
of \(n_\omega\) (along with \(R_{\mathrm{joint}}L_R\) and the residual suprema
over \(\mathcal F_{\mathrm{joint}}\)). Consequently, when \(M\) is large
relative to \(n_\omega\) and the corresponding radii, Lipschitz constants, and
residual suprema are comparable, the sampling error bound for
\(\mathcal F_\star\) can be substantially smaller than that for
\(\mathcal F_{\mathrm{joint}}\).

By decoupling coefficient fitting from parameter selection, the bilevel
formulation uses the validation points only to select the \(n_\omega\)
parameters in \(\boldsymbol\omega\), which yields a tighter validation sampling bound than joint fitting on
the same validation data. As shown by the numerical comparison in Example~A of Section \ref{sec:numerical-experiments},
bilevel selection produces a more accurate PDE solution than single-level
joint optimization.

\subsection{Validation and Residual-Stability Verification for Nonlinear Elliptic Kansa Solvers}
\label{subsec:exampleA-verification}

The preceding analysis requires two problem-dependent ingredients: a uniform
bound on the validation sampling error over the local residual class \(\mathcal F_{\mathrm{loc}}\)
and control of the PDE solution error by the corresponding population residual
norm. We verify these requirements for the unregularized Gaussian Kansa
discretization of \(-\Delta u+u^3=f\) on \(\Omega=(0,1)^2\) with homogeneous Dirichlet boundary
conditions in Example~A of Section \ref{sec:numerical-experiments}, thereby obtaining an \(L^2(\Omega)\) solution error
bound from Corollary~\ref{cor:local-class-a-posteriori}.

We use the scalar log-lengthscale \(\omega\), with \(\ell=e^\omega\).
Before drawing the validation set, fix \(R_\beta>0\),
\(-\infty<\omega_-<\omega_+<\infty\), and a nonempty compact set
\(\mathcal K_A\) contained in
\[
\left\{(\boldsymbol\beta,\omega):
\|\boldsymbol\beta\|\le R_\beta,
\quad
\omega_-\le\omega\le\omega_+\right\}.
\]
Assume that \(f\in C(\overline\Omega)\). Since the Gaussian basis
functions and their spatial derivatives depend smoothly on \(\omega\) and the
trial function is linear in \(\boldsymbol\beta\), the residual map is
continuous from \(\mathcal K_A\) into
\(C(\overline\Omega)\times C(\partial\Omega)\). Hence
Proposition~\ref{prop:local-class-validation-sampling} applies with
\(\mathcal K_{\mathrm{loc}}=\mathcal K_A\).

The validation points in Example~A are sampled uniformly, so
\(d\nu_I=|\Omega|^{-1}\,d\boldsymbol x\) and
\(d\nu_B=|\partial\Omega|^{-1}\,ds\).
Take \(X=L^2(\Omega)\), \(Y_I=L^2(\Omega)\), and \(Y_B=L^2(\partial\Omega)\).
Since
\(\|r\|_Y^2=|\Omega|\,\|r_I\|_{L^2(\nu_I)}^2
+\eta_B|\partial\Omega|\,\|r_B\|_{L^2(\nu_B)}^2\),
every residual pair \(r=(r_I,r_B)\) satisfies
\begin{equation}
\label{eq:exampleA-local-norm-comparison}
\|r\|_Y
\le
C_{\mathrm{loc}}\|r\|_{Y_0},
\qquad
C_{\mathrm{loc}}^2
=
\max\left\{|\Omega|,\eta_B|\partial\Omega|\right\},
\end{equation}
so \eqref{eq:residual-norm-compatibility} holds. It remains to establish the residual-stability estimate
\eqref{eq:continuous-residual-stability}.

\begin{proposition}
\label{prop:nonlinear-elliptic-residual-stability}
Let \(\Omega=(0,1)^2\), let \(f\in L^2(\Omega)\), and let
\(u\in H_0^1(\Omega)\) be the unique weak solution of
\[
-\Delta u+u^3=f
\quad\text{in }H^{-1}(\Omega).
\]
For every \(R>0\), there is \(C_R>0\) such that each
\(v\in H^1(\Omega)\) satisfying \(\|v\|_{H^1(\Omega)}\le R\) and
\[
r_I(v):=-\Delta v+v^3-f\in L^2(\Omega),
\qquad
r_B(v):=v|_{\partial\Omega},
\]
obeys
\begin{equation*}
\|v-u\|_{L^2(\Omega)}
\le
C_R\left(
\|r_I(v)\|_{L^2(\Omega)}
+
\|r_B(v)\|_{L^2(\partial\Omega)}
\right).
\end{equation*}
\end{proposition}

\begin{proof}
Set \(e=v-u\) and
\(q=v^2+vu+u^2\). Then \(q\ge0\) and
\[
-\Delta e+qe=r_I(v),
\qquad e|_{\partial\Omega}=r_B(v).
\]
The two-dimensional Sobolev embeddings imply
\(q\in L^4(\Omega)\) and \(qe=v^3-u^3\in L^2(\Omega)\), with
\(\|q\|_{L^4(\Omega)}\) bounded in terms of \(R\) and
\(\|u\|_{H^1(\Omega)}\); hence \(\Delta e\in L^2(\Omega)\).
To bound \(\|e\|_{L^2(\Omega)}\) in terms of the interior and boundary
residuals, let \(z\in H_0^1(\Omega)\) solve the dual problem
\[
-\Delta z+qz=e.
\]
Since \(q\ge0\), coercivity and the embedding
\(H_0^1(\Omega)\hookrightarrow L^4(\Omega)\), together with the bound on \(\|q\|_{L^4(\Omega)}\),
control \(\|z\|_{H^1(\Omega)}\) and \(\|qz\|_{L^2(\Omega)}\) by
\(\|e\|_{L^2(\Omega)}\), so that \(-\Delta z=e-qz\in L^2(\Omega)\).
Dirichlet \(H^2\)-regularity on the convex polygon \(\Omega\) and the normal-trace
estimate then yield
\begin{equation}
\label{eq:dual-regularity-compressed}
\|z\|_{H^2(\Omega)}
+\|\partial_n z\|_{L^2(\partial\Omega)}
\le C_R\|e\|_{L^2(\Omega)}.
\end{equation}
Applying Green's identity with \(z|_{\partial\Omega}=0\) and using the dual equation gives
\[
\begin{aligned}
\|e\|_{L^2(\Omega)}^2
&=\int_\Omega(-\Delta e+qe)z\,d\boldsymbol x
-\int_{\partial\Omega}r_B(v)\,\partial_n z\,ds\\
&=\int_\Omega r_I(v)z\,d\boldsymbol x
-\int_{\partial\Omega}r_B(v)\,\partial_n z\,ds.
\end{aligned}
\]
Consequently, the Cauchy--Schwarz inequality and \eqref{eq:dual-regularity-compressed} give
\[
\|e\|_{L^2(\Omega)}^2
\le
C_R\bigl(\|r_I(v)\|_{L^2(\Omega)}
+\|r_B(v)\|_{L^2(\partial\Omega)}\bigr)
\|e\|_{L^2(\Omega)}.
\]
Division by \(\|e\|_{L^2(\Omega)}\) proves the claim when \(e\ne0\); the
case \(e=0\) is immediate.
\end{proof}

To apply Proposition~\ref{prop:nonlinear-elliptic-residual-stability} on
\(\mathcal K_A\), note that the parameter bounds on \(\mathcal K_A\) yield a
constant \(R_A>0\) such that
\(\|\widehat u(\cdot;\boldsymbol\beta,\omega)\|_{H^1(\Omega)}\le R_A\) for all
\((\boldsymbol\beta,\omega)\in\mathcal K_A\).
Since \(\|r_I\|_{L^2(\Omega)}+\|r_B\|_{L^2(\partial\Omega)}
\le(1+\eta_B^{-1})^{1/2}\|r\|_Y\) by
\eqref{eq:continuous-residual-norm},
Proposition~\ref{prop:nonlinear-elliptic-residual-stability} with \(R=R_A\)
verifies the residual-stability estimate
\eqref{eq:continuous-residual-stability} on \(\mathcal K_A\) with
\(C_{\mathrm{stab}}=C_{R_A}(1+\eta_B^{-1})^{1/2}\).
Together with \eqref{eq:exampleA-local-norm-comparison}, this verifies
Assumption~\ref{ass:residual-stability-compatibility} on \(\mathcal K_A\).
Corollary~\ref{cor:local-class-a-posteriori} therefore gives, on the uniform
sampling event~\eqref{eq:local-class-validation-sampling},
\[
\|\widehat{u}(\cdot;\boldsymbol\beta,\omega)-u\|_{L^2(\Omega)}
\le
C_{\mathrm{stab}}C_{\mathrm{loc}}
\left(
2\mathcal J(\boldsymbol\beta,\omega;\mathcal V)
+
\delta_{\mathrm{samp}}^2
\right)^{1/2},
\qquad
(\boldsymbol\beta,\omega)\in\mathcal K_A.
\]

\subsection{Terminal and Finite-Iteration PDE Error Bounds}
\label{subsec:end-to-end-transfer}

Theorem~\ref{thm:bilevel-validation-oracle} bounds the PDE error at a
validation-selected parameter value when the selected lower-level
stationary response is evaluated exactly, whereas in practice the
lower-level solve returns an approximate coefficient vector. To convert the resulting coefficient error into a solution error in
\(X\), we define, for \(\boldsymbol\omega\in\mathcal W\) and
\(\boldsymbol a\in\mathbb R^M\),
\[
\Psi_{\boldsymbol\omega}\boldsymbol a
:=
\sum_{j=1}^M a_j\psi_j(\cdot;\boldsymbol\omega),
\qquad
L_{\Psi,X}
:=
\sup_{\boldsymbol\omega\in\mathcal W,\ \|\boldsymbol a\|=1}
\|\Psi_{\boldsymbol\omega}\boldsymbol a\|_X.
\]
Recall from Section~\ref{subsec:validation-residual-error-control} that
\(R_\star=\inf_{\boldsymbol\omega\in\mathcal W}
\|\mathcal R(\widehat{u}(\cdot;\boldsymbol\beta^\star(\boldsymbol\omega),
\boldsymbol\omega))\|_{Y_0}^2\)
is the best squared population residual norm along the selected response
branch. Proposition~\ref{prop:sc-local-response-estimates} bounds the coefficient
error in terms of the lower-level stationarity residual. Combining this estimate
with the definition of \(L_{\Psi,X}\) and
Theorem~\ref{thm:bilevel-validation-oracle} gives the following terminal PDE
error bound.

\begin{theorem}
\label{thm:end-to-end-terminal-transfer}
Suppose the hypotheses of Theorem~\ref{thm:bilevel-validation-oracle}
hold for \(\widetilde{\boldsymbol\omega}\), and \(L_{\Psi,X}<\infty\).
Assume that Proposition~\ref{prop:sc-local-response-estimates} applies at
\((\boldsymbol\beta^\star(\widetilde{\boldsymbol\omega}),
\widetilde{\boldsymbol\omega})\), with neighborhood
\(\mathcal N_{\widetilde{\boldsymbol\omega}}\) and constant \(m_{\widetilde{\boldsymbol\omega}}>0\).
Let \(\widehat{\boldsymbol\beta}\in\mathcal N_{\widetilde{\boldsymbol\omega}}\) satisfy
\begin{equation}
\label{eq:terminal-lower-stationarity}
\|\nabla_{\boldsymbol\beta}\mathcal E(
\widehat{\boldsymbol\beta},\widetilde{\boldsymbol\omega})\|
\le \varepsilon_\beta.
\end{equation}
On the event \eqref{eq:bilevel-response-class-validation-sampling},
\begin{equation}
\label{eq:terminal-pde-transfer}
\begin{aligned}
\|\widehat{u}(\cdot;\widehat{\boldsymbol\beta},
\widetilde{\boldsymbol\omega})-u\|_X
\le{}&
C_{\mathrm{stab}}C_{\mathrm{loc}}
\left[
R_\star+2\delta_{\mathrm{samp},\star}^2
+2\varepsilon_{\mathrm{out}}
\right]^{1/2}
+
\frac{L_{\Psi,X}}{m_{\widetilde{\boldsymbol\omega}}}\varepsilon_\beta.
\end{aligned}
\end{equation}
\end{theorem}

\begin{proof}
Theorem~\ref{thm:bilevel-validation-oracle} bounds the exact-response error by
the first term on the right-hand side of
\eqref{eq:terminal-pde-transfer}. Proposition~\ref{prop:sc-local-response-estimates}
and \eqref{eq:terminal-lower-stationarity} give
\[
\|\widehat{\boldsymbol\beta}
-\boldsymbol\beta^\star(\widetilde{\boldsymbol\omega})\|
\le
\frac{\varepsilon_\beta}{m_{\widetilde{\boldsymbol\omega}}}.
\]
Linearity of the trial representation therefore yields
\[
\|\widehat{u}(\cdot;\widehat{\boldsymbol\beta},\widetilde{\boldsymbol\omega})
-\widehat{u}(\cdot;\boldsymbol\beta^\star(\widetilde{\boldsymbol\omega}),
\widetilde{\boldsymbol\omega})\|_X
\le
\frac{L_{\Psi,X}}{m_{\widetilde{\boldsymbol\omega}}}\varepsilon_\beta.
\]
The triangle inequality proves \eqref{eq:terminal-pde-transfer}.
\end{proof}

The bound \eqref{eq:terminal-pde-transfer} separates the effect of the
upper-level objective gap \(\varepsilon_{\mathrm{out}}\) from that of the
terminal lower-level stationarity residual \(\varepsilon_\beta\). Because the
upper-level optimization enters \eqref{eq:terminal-pde-transfer} only through
\(\varepsilon_{\mathrm{out}}\), Theorem~\ref{thm:end-to-end-terminal-transfer}
applies to any upper-level method, such as Adam or AdamW in
Section~\ref{sec:numerical-experiments}.

For the basic
upper-level update \eqref{eq:npde:gn-update}, \(\varepsilon_{\mathrm{out}}\) can be bounded
explicitly. We apply Theorem~\ref{thm:end-to-end-terminal-transfer} at
\(\widetilde{\boldsymbol\omega}=\boldsymbol\omega^K\) after \(K\) iterations,
taking \(\varepsilon_{\mathrm{out}}=\Delta_K\), where the upper-level
objective gap and hypergradient error at iteration \(k\) are
\[
\Delta_k
:=
\Phi(\boldsymbol\omega^k)-\Phi_\star,
\qquad
\boldsymbol e^k
:=
\boldsymbol d_s^k-\nabla\Phi(\boldsymbol\omega^k).
\]
Suppose that \(\Phi\) satisfies the local Polyak--\L{}ojasiewicz (PL) condition with constant
\(\mu_\Phi>0\):
\begin{equation}
\label{eq:end-to-end-pl}
\frac{1}{2}\|\nabla\Phi(\boldsymbol\omega)\|^2
\ge
\mu_\Phi
\bigl(\Phi(\boldsymbol\omega)-\Phi_\star\bigr),
\qquad
\boldsymbol\omega\in\mathcal W.
\end{equation}
For \(0\le a\le b\le K\), define the cumulative contraction factors
\[
Q_{a,b}
:=
\prod_{i=a}^{b-1}(1-\alpha_i\mu_\Phi),
\qquad
Q_{b,b}:=1.
\]
An inexact-gradient descent recursion bounds the terminal gap
\(\Delta_K\) by the contracted initial gap and the accumulated hypergradient
errors. Combining this estimate with Theorem~\ref{thm:end-to-end-terminal-transfer}
yields the following finite-iteration PDE error bound.

\begin{theorem}
\label{thm:end-to-end-pde-error}
Suppose the hypotheses of
Corollary~\ref{cor:bilevel-response-validation-sampling} and
Assumption~\ref{ass:residual-stability-compatibility} hold,
\(\mathcal G_\star\subset\mathcal K_{\mathrm{loc}}\),
and \(L_{\Psi,X}<\infty\).
Let \(\{\boldsymbol\omega^k\}_{k=0}^{K}\subset\mathcal W\) be generated by
\eqref{eq:npde:gn-update}, with step sizes satisfying
\[
0<\alpha_k\le L_\Phi^{-1},
\qquad
0<\alpha_k\mu_\Phi<1.
\]
Assume that every update segment lies in a region where \(\nabla\Phi\) is
\(L_\Phi\)-Lipschitz, and that \eqref{eq:end-to-end-pl} holds on \(\mathcal W\).
At the terminal iterate \(\boldsymbol\omega^K\), assume the local lower-level
response conditions of Theorem~\ref{thm:end-to-end-terminal-transfer} hold
for \(\widehat{\boldsymbol\beta}^K\), with constant \(m_K>0\) and
\(\|\nabla_{\boldsymbol\beta}\mathcal E(
\widehat{\boldsymbol\beta}^K,\boldsymbol\omega^K)\|\le\varepsilon_K\).
Then, on the event
\eqref{eq:bilevel-response-class-validation-sampling},
\begin{equation}
\label{eq:end-to-end-finite-pde-error}
\begin{aligned}
\|\widehat{u}(\cdot;\widehat{\boldsymbol\beta}^K,
\boldsymbol\omega^K)-u\|_X
\le{}
C_{\mathrm{stab}}C_{\mathrm{loc}}
\Bigg[
R_\star
+2\delta_{\mathrm{samp},\star}^2
+2Q_{0,K}\Delta_0
+\sum_{j=0}^{K-1}
\alpha_jQ_{j+1,K}\|\boldsymbol e^j\|^2
\Bigg]^{1/2} +
\frac{L_{\Psi,X}}{m_K}\varepsilon_K.
\end{aligned}
\end{equation}
\end{theorem}

\begin{proof}
The descent
lemma, \(L_\Phi\alpha_k\le1\), and
\(\boldsymbol e^k=\boldsymbol d_s^k-
\nabla\Phi(\boldsymbol\omega^k)\) give
\[
\Phi(\boldsymbol\omega^{k+1})
\le
\Phi(\boldsymbol\omega^k)
-\frac{\alpha_k}{2}\|\nabla\Phi(\boldsymbol\omega^k)\|^2
+\frac{\alpha_k}{2}\|\boldsymbol e^k\|^2.
\]
Subtracting \(\Phi_\star\), applying \eqref{eq:end-to-end-pl}, and iterating
yield
\begin{equation}
\label{eq:end-to-end-outer-gap}
\Delta_K
\le
Q_{0,K}\Delta_0
+\frac{1}{2}\sum_{j=0}^{K-1}
\alpha_jQ_{j+1,K}\|\boldsymbol e^j\|^2.
\end{equation}
Thus \(\boldsymbol\omega^K\) satisfies
\eqref{eq:bilevel-approx-upper-optimality} with
\(\varepsilon_{\mathrm{out}}=\Delta_K\). Apply
Theorem~\ref{thm:end-to-end-terminal-transfer} and use
\eqref{eq:end-to-end-outer-gap}.
\end{proof}

In \eqref{eq:end-to-end-finite-pde-error}, each squared
hypergradient error \(\|\boldsymbol e^j\|^2\)
is weighted by
\(\alpha_jQ_{j+1,K}\), so the contributions from early steps decay across
subsequent iterations. Under the hypotheses of
Theorem~\ref{thm:sc-hypergradient-decomposition},
\eqref{eq:sc:refined-main-bias} bounds each \(\|\boldsymbol e^j\|\)
in terms of the lower-level and final adjoint residuals and the damping \(\mu\).
For a constant step size, the initial-gap term in \eqref{eq:end-to-end-finite-pde-error}
vanishes as \(K\to\infty\), and when these residuals meet fixed tolerances, passing to the limit gives the following
asymptotic PDE error bound.

\begin{corollary}
\label{cor:end-to-end-asymptotic-floor}
Suppose the hypotheses of Theorem~\ref{thm:end-to-end-pde-error} hold
for every \(K\ge1\) with a constant step size \(\alpha_k\equiv\alpha\).
Assume further that the
hypotheses of Theorem~\ref{thm:sc-hypergradient-decomposition} hold for every
\(k\ge0\), that the lower-level and final adjoint residuals satisfy
\[
\eta_k=\|\nabla_{\boldsymbol\beta}\mathcal E(
\widehat{\boldsymbol\beta}^k,\boldsymbol\omega^k)\|
\le\varepsilon_\beta,
\qquad
\|\mathbf r_s^k\|\le\tau_{\rm adj}\max\{1,\|\boldsymbol b^k\|\}
\qquad(k\ge0)
\]
for some \(\varepsilon_\beta\ge0\), and that
\(C_{3,k}\le\overline C_3\) and \(m_k\ge\underline m\) in
\eqref{eq:sc:refined-bias-via-tolerances} for all \(k\ge0\), with
\(\overline C_3,\underline m>0\).
Then, on the event
\eqref{eq:bilevel-response-class-validation-sampling},
\begin{equation}
\label{eq:asymptotic-pde-bound}
\begin{aligned}
\limsup_{k\to\infty}
\|\widehat{u}(\cdot;\widehat{\boldsymbol\beta}^k,
\boldsymbol\omega^k)-u\|_X
\le{}&
C_{\mathrm{stab}}C_{\mathrm{loc}}
\left[
R_\star
+2\delta_{\mathrm{samp},\star}^2
+\frac{\overline C_3^2}{\mu_\Phi}
\left(
\frac{\varepsilon_\beta}{\underline m}
+\mu+\tau_{\rm adj}
\right)^2
\right]^{1/2}
+
\frac{L_{\Psi,X}}{\underline m}
\varepsilon_\beta.
\end{aligned}
\end{equation}
\end{corollary}

\begin{proof}
By \eqref{eq:sc:refined-bias-via-tolerances} with
\(\varepsilon_k=\varepsilon_\beta\), together with
\(C_{3,k}\le\overline C_3\) and \(m_k\ge\underline m\), every iterate
\(j\ge0\) satisfies
\[
\|\boldsymbol e^j\|
\le
\overline C_3\left(
\frac{\varepsilon_\beta}{\underline m}+\mu+\tau_{\rm adj}
\right).
\]
For \(\alpha_j\equiv\alpha\) with \(0<\alpha\mu_\Phi<1\), we have
\(Q_{j+1,k}=(1-\alpha\mu_\Phi)^{k-1-j}\), and hence, for every \(k\ge1\),
\[
\sum_{j=0}^{k-1}
\alpha Q_{j+1,k}\|\boldsymbol e^j\|^2
\le
\frac{\overline C_3^2}{\mu_\Phi}
\left(
\frac{\varepsilon_\beta}{\underline m}+\mu+\tau_{\rm adj}
\right)^2
\left[1-(1-\alpha\mu_\Phi)^k\right].
\]
Apply \eqref{eq:end-to-end-finite-pde-error} at iteration \(k\) with
\(\varepsilon_k=\varepsilon_\beta\), use \(m_k\ge\underline m\) and
\(Q_{0,k}=(1-\alpha\mu_\Phi)^k\), and take the limit superior as
\(k\to\infty\).
\end{proof}

When
\(\varepsilon_\beta=\mu=\tau_{\rm adj}=0\), the asymptotic PDE error bound
\eqref{eq:asymptotic-pde-bound} reduces
to the floor
\(C_{\mathrm{stab}}C_{\mathrm{loc}}(R_\star+2\delta_{\mathrm{samp},\star}^2)^{1/2}\),
which is the bound \eqref{eq:bilevel-pde-oracle} of
Theorem~\ref{thm:bilevel-validation-oracle} with \(\varepsilon_{\mathrm{out}}=0\).
With inexact solves, \eqref{eq:asymptotic-pde-bound} includes both the
hypergradient-error term
\(\tfrac{\overline C_3^2}{\mu_\Phi}(\varepsilon_\beta/\underline m+\mu+\tau_{\rm adj})^2\)
inside the square root and the terminal fitting term
\(\tfrac{L_{\Psi,X}}{\underline m}\varepsilon_\beta\), which reflects the
accuracy of the final lower-level solve. Even when the terminal fitting term is small, the
hypergradient-error term can keep the bound in
\eqref{eq:asymptotic-pde-bound} well above the floor, so controlling the
hypergradient error remains essential.

Under the assumptions of
Corollary~\ref{cor:end-to-end-asymptotic-floor}, the estimates in
Section~\ref{sec:bias} provide explicit control of both contributions through
the tolerances \((\varepsilon_\beta,\tau_{\rm adj})\) and the damping \(\mu\)
used in Algorithm~\ref{alg:npde-gn}. Once \((\varepsilon_\beta,\mu,\tau_{\rm adj})\)
make the hypergradient-error term small relative to
\(R_\star+2\delta_{\mathrm{samp},\star}^2\) and the terminal fitting
term small relative to the
floor, solving the lower-level and adjoint problems to higher accuracy yields
little further reduction in the PDE error bound, which explains why inexact
solves suffice.

\section{Applications and Numerical Experiments}
\label{sec:numerical-experiments}

This section presents five nonlinear PDE examples to evaluate
Algorithm \ref{alg:npde-gn} across Kansa, RFM/ELM-type,
and PINN-type mesh-free solvers. Unless otherwise stated, all numerical
experiments were
conducted using Python~3.10.19 on a Mac Studio equipped with
an Apple M4 Max processor and running macOS.
The code for the numerical experiments is available at
\url{https://github.com/SUSTech-Optimization/Bihyperparameter_PDE}.

The reported training and validation residuals are
\(\|\boldsymbol\rho_{\mathcal T}^{\rm tr}(\boldsymbol\beta,\boldsymbol\omega)\|_2\)
and \(\|\boldsymbol\rho(\boldsymbol\beta,\boldsymbol\omega;\mathcal V)\|_2\),
respectively.
The relative \(L^2\) and absolute \(L^\infty\) errors are computed at
the evaluation points \(\{\boldsymbol x_q\}_{q=1}^{N_{\mathrm{eval}}}\)
specified in each example, with \(u_{\mathrm{ref}}\) denoting the exact or
reference solution:
\[
\operatorname{RelErr}_{L^2}
:=
\frac{
\left(\sum_{q=1}^{N_{\mathrm{eval}}}
|\widehat{u}(\boldsymbol x_q)-u_{\mathrm{ref}}(\boldsymbol x_q)|^2\right)^{1/2}
}{
\left(\sum_{q=1}^{N_{\mathrm{eval}}}
|u_{\mathrm{ref}}(\boldsymbol x_q)|^2\right)^{1/2}
},
\qquad
\operatorname{Err}_{L^\infty}
:=
\max_{1\le q\le N_{\mathrm{eval}}}
|\widehat{u}(\boldsymbol x_q)-u_{\mathrm{ref}}(\boldsymbol x_q)|.
\]

For the lower-level solves in Algorithm~\ref{alg:npde-gn}, we use SciPy's
\texttt{least\_squares} in Examples~A--C and~E. In Example~A, the budget is
\(T=100\) residual evaluations for the first solve and \(T=50\) for each
subsequent solve. In Examples~B, C, and~E, the budget is \(T=100\) throughout.
In Example~D, we use a JAXopt solver with \(T=10\). The lower-level stopping tolerance is set to
\(\varepsilon_k=10^{-3}\) in Example~A and \(\varepsilon_k=10^{-6}\) in
Examples~B--E.

To accelerate the upper-level optimization, we
replace the basic update \eqref{eq:npde:gn-update} by a full-batch
Adam or AdamW update computed from the approximate hypergradient direction returned by
Algorithm~\ref{alg:npde-implicit-hypergradient}, with the choice of optimizer and the
adjoint parameters \((\mu,\tau_{\mathrm{adj}},s)\) specified in each example. The pointwise
hypergradient analysis in Section~\ref{sec:bias} and the terminal PDE error
bound in Theorem~\ref{thm:end-to-end-terminal-transfer} do not require a
specific upper-level optimizer; the finite-iteration analysis in
Theorem~\ref{thm:end-to-end-pde-error} is developed for the basic
update.

\subsection{Example A (Kansa): Nonlinear elliptic problem}
\label{subsec:exampleA}

We consider the nonlinear elliptic boundary-value problem on \(\Omega=(0,1)^2\),
\[
-\Delta u(\boldsymbol{x}) + u(\boldsymbol{x})^3 = f(\boldsymbol{x}),
\qquad \boldsymbol{x}\in\Omega,
\qquad
u(\boldsymbol{x}) = 0,
\qquad \boldsymbol{x}\in\partial\Omega.
\]
The exact solution is
\[
u(x,y)=\sin(\pi x)\sin(\pi y)+4\sin(4\pi x)\sin(4\pi y),
\]
and the forcing term is chosen accordingly as
$
f(x,y)=-\Delta u(x,y)+\bigl(u(x,y)\bigr)^3.
$

The Kansa approximation uses Gaussian kernels with a shared lengthscale
\(\ell=e^\omega\):
\[
\widehat{u}(\boldsymbol{x};\boldsymbol{\beta},\omega)
=
\sum_{j=1}^{M}\beta_j
\exp\!\left(
-\frac{\|\boldsymbol{x}-\boldsymbol{c}_j\|_2^2}{\ell^2}
\right).
\]
We use the bilevel formulation
\eqref{eq:outer_objective}--\eqref{eq:inner_constraint} with \(\lambda=0\)
and select the scalar log-lengthscale \(\omega\).

\begin{figure}[!htb]
    \centering
    \subcaptionbox{Exact solution\label{fig:exampleA-solution-exact}}[0.32\textwidth]{%
        \includegraphics[width=\linewidth]{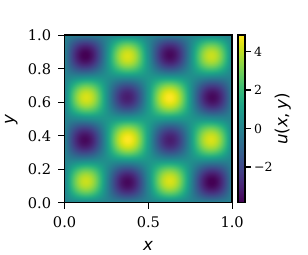}}
    \hfill
    \subcaptionbox{Kansa solution\label{fig:exampleA-solution-kansa}}[0.32\textwidth]{%
        \includegraphics[width=\linewidth]{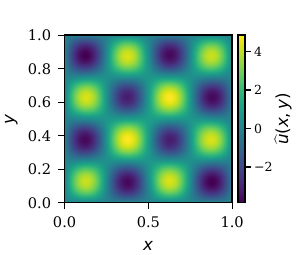}}
    \hfill
    \subcaptionbox{Absolute error\label{fig:exampleA-solution-error}}[0.32\textwidth]{%
        \includegraphics[width=\linewidth]{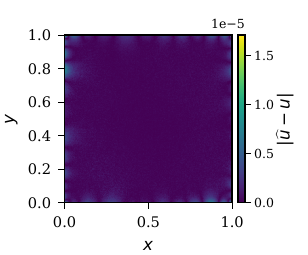}}
    \caption{Kansa solution after \(100\)
    bilevel steps for Example A.}
    \label{fig:exampleA-solution}
\end{figure}

The experiment uses \(640\) Kansa centers, including \(512\) interior and
\(128\) boundary centers. The training and validation sets each contain
\(2048\) interior and \(512\) boundary collocation points sampled uniformly,
and errors are evaluated on a uniform \(256\times256\) spatial
grid. Starting from
\(\ell_0=0.1\), we perform \(100\) bilevel steps with Adam and
step size \(0.1\). The adjoint computation uses
\(\mu=10^{-6}\),
\(\tau_{\mathrm{adj}}=10^{-6}\), and \(s=10\).
Figure~\ref{fig:exampleA-solution} shows the Kansa solution
after \(100\) bilevel steps and the pointwise absolute
error, which remains below \(2\times10^{-5}\).

To evaluate the proposed method, we compare it with two
baselines that use the same initial lengthscale, Kansa centers, and
training/validation sets.
The first baseline follows the GN-linearization method in~\cite{nelsen2025bilevel},
which linearizes the training and validation residuals about the current Kansa
solution at each step. We substitute the resulting closed-form coefficient
response into the linearized validation objective and take one Adam step in
\(\omega\) on the reduced objective. After \(1500\) such steps, we refit the
coefficients at the selected lengthscale using the original nonlinear training
objective, starting from the final linearized coefficient vector.
The second is single-level joint optimization, which minimizes the lower-level training
objective jointly over \((\boldsymbol\beta,\omega)\) using the same
GN-type solver, with a budget of \(3000\).

\begin{figure}[!htb]
    \centering
    \subcaptionbox{Hypergradient comparison\label{fig:exampleA-hypergradient-panel}}[0.48\textwidth]{%
        \includegraphics[width=\linewidth]{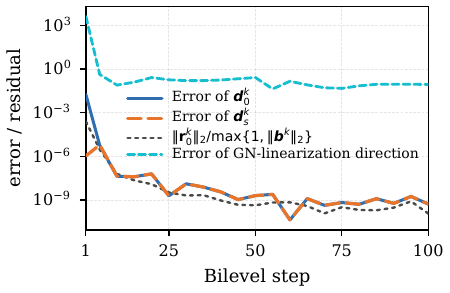}}
    \hfill
    \subcaptionbox{Errors and residuals of the proposed method (\(s=10\))\label{fig:exampleA-residual-panel}}[0.48\textwidth]{%
        \includegraphics[width=\linewidth]{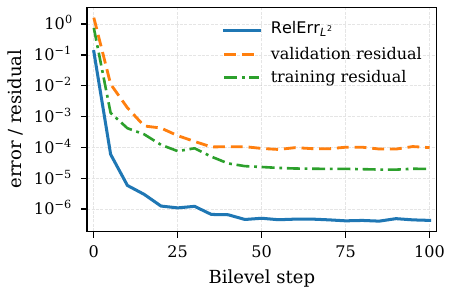}}
    \par\medskip
    \subcaptionbox{Error vs.\ wall-clock time\label{fig:exampleA-lengthscale-panel}}[0.48\textwidth]{%
        \includegraphics[width=\linewidth]{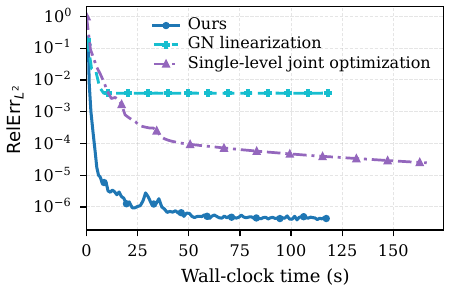}}
    \hfill
    \subcaptionbox{Lengthscale study\label{fig:exampleA-sweep-panel}}[0.48\textwidth]{%
        \includegraphics[width=\linewidth]{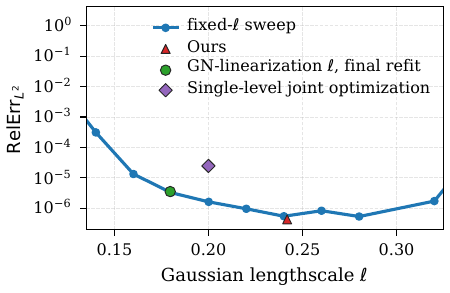}}
    \caption{Hypergradient, accuracy, and parameter-selection
    diagnostics for Example A.}
    \label{fig:exampleA-outer-summary}
\end{figure}

We first examine the accuracy of the hypergradient directions.
To isolate
the error from the adjoint computation, we compare three directions with damped
exact-Hessian references at their respective computed responses
in Figure~\ref{fig:exampleA-outer-summary}(a): \(\boldsymbol d_0^k\) and
\(\boldsymbol d_s^k\), computed from \eqref{eq:npde:gn-hg-refined} using the
initial and returned adjoints, respectively, and the GN-linearization direction. The
reference adjoint systems are solved to high accuracy by a double-precision
Bunch--Kaufman factorization with iterative refinement.
The error of \(\boldsymbol d_0^k\)
closely follows the normalized initial adjoint residual
\(\|\mathbf r_0^k\|_2/\max\{1,\|\boldsymbol b^k\|_2\}\) along the trajectory, supporting the use of this residual as an error indicator
for \(\boldsymbol d_0^k\).
The residual exceeds \(\tau_{\mathrm{adj}}\) only at the first few
steps, where the MINRES iterations substantially reduce these
larger errors; throughout the trajectory, \(\boldsymbol d_s^k\)
has an error several orders of magnitude smaller than that of the
GN-linearization direction.

We next examine the relation between the validation
residual and the PDE solution error along the iterations of the proposed method.
As shown in Figure~\ref{fig:exampleA-outer-summary}(b),
\(\operatorname{RelErr}_{L^2}\)
decreases alongside the training and validation residuals, which is
consistent with the residual-to-solution analysis in
Section~\ref{sec:validation-residual-analysis} and supports the use of the
validation residual as a selection criterion for this problem.

Finally, we assess the computational efficiency and parameter-selection
quality of the three methods. Figure~\ref{fig:exampleA-outer-summary}(c) shows
that the proposed method reduces \(\operatorname{RelErr}_{L^2}\) below
\(10^{-5}\) within about \(10\)~s, whereas the GN-linearization error remains
near \(4\times10^{-3}\) before the final refit, and single-level joint
optimization reaches \(2.5\times10^{-5}\) only after about \(165\)~s.
Table~\ref{tab:exampleA-lengthscale} and
Figure~\ref{fig:exampleA-outer-summary}(d) compare the selected lengthscales
and \(\operatorname{RelErr}_{L^2}\) of the three methods
against fixed-lengthscale solves. The proposed method selects \(\ell=0.2417\), within the low-error range of the fixed-lengthscale results,
and attains the smallest error in the table. GN linearization selects a
smaller lengthscale, \(\ell=0.1795\), and its refitted error is nearly an
order of magnitude larger than that of the proposed method. Single-level joint
optimization yields the largest error.

\begin{table}[!htb]
\centering
\small
\begin{tabular}{lcc}
\hline
Method & Lengthscale & $\operatorname{RelErr}_{L^2}$ \\
\hline
Fixed & \(\ell=0.16\) & \(1.3274\times10^{-5}\) \\
Fixed & \(\ell=0.20\) & \(1.6174\times10^{-6}\) \\
Fixed & \(\ell=0.24\) & \(5.4363\times10^{-7}\) \\
Fixed & \(\ell=0.28\) & \(5.3204\times10^{-7}\) \\
Fixed & \(\ell=0.32\) & \(1.7161\times10^{-6}\) \\
\hline
GN-linearization \(\ell\), final refit & \(\ell=0.1795\) & \(3.5236\times10^{-6}\) \\
Single-level joint optimization & \(\ell=0.1998\) & \(2.4566\times10^{-5}\) \\
Ours, \(s=10\) & \(\ell=0.2417\) & \(\mathbf{4.3057\times10^{-7}}\) \\
\hline
\end{tabular}
\caption{Lengthscales and \(\operatorname{RelErr}_{L^2}\) at fixed lengthscales and for the compared methods in Example A.}\label{tab:exampleA-lengthscale}
\end{table}

\FloatBarrier 

\subsection{Example B (RFM/ELM): Nonlinear Helmholtz problem}
\label{subsec:exampleB}

Following the nonlinear Helmholtz benchmark used in SO-PIFRNN \cite{linghu2025so}, we consider the problem on \(\Omega=(0,1)^2\),
\begin{equation}\label{eq:helmholtz}
    \Delta u(\boldsymbol{x})-100u(\boldsymbol{x})+10\cos(u(\boldsymbol{x}))=f(\boldsymbol{x}),
\qquad
\boldsymbol{x}\in\Omega.
\end{equation}
The exact solution is
\[
u(x,y)=4\cos(3\pi x^2)\sin(3\pi y^2),
\]
and the Dirichlet data \(g=u|_{\partial\Omega}\) and the forcing term
\(f(x,y)=\Delta u(x,y)-100u(x,y)+10\cos\bigl(u(x,y)\bigr)\) are chosen
accordingly.

We use the sine random features
\(\psi_j(\boldsymbol x;R_m)
=\sin\bigl(R_m(\widetilde{\boldsymbol w}_j\cdot\boldsymbol x+\widetilde b_j)\bigr)\),
as in standard RFM/ELM constructions \cite{huang2006extreme,dong2021local}, and
sample the base variables as
\[
\widetilde{\boldsymbol w}_j\sim\mathcal U([-1,1]^2),
\qquad
\widetilde b_j\sim\mathcal U([-1,1]),
\qquad j=1,\ldots,M.
\]

We use the bilevel formulation
\eqref{eq:outer_objective}--\eqref{eq:inner_constraint} with
\(\mathbf W_{\mathcal T}(\boldsymbol\omega)
=\operatorname{diag}(\mathbf I,\alpha_B\mathbf I)\) and \(\mathbf R_{\mathcal T}=\mathbf I\).

The solver-design parameters \(R_m\), \(\alpha_B\), and \(\lambda\) determine
the feature scale, the training boundary-residual weight,
and the ridge-regularization strength, respectively. To preserve positivity, we
optimize them in log coordinates,
\(\boldsymbol\omega
=(\log R_m,\log\alpha_B,\log\lambda)\).

The experiment uses \(M=1000\) sine random features. The
training and validation sets each contain
\(1000\) interior and \(1000\) boundary points, and errors are evaluated on a
uniform \(256\times256\) spatial grid. Starting from
\((R_m,\alpha_B,\lambda)=(20,1,10^{-3})\), we run \(40\) bilevel steps with Adam, using
initial step size \(0.5\) and multiplicative decay \(0.95\). The
adjoint computation uses \(\mu=0\),
\(\tau_{\mathrm{adj}}=10^{-8}\), and \(s=5\).

To evaluate the proposed method, we compare it with three baselines that use
the same initial parameters, training/validation sets, and number of random
features. As in Example~A, the GN-linearization baseline linearizes the
training and validation residuals at the current solution and takes one Adam
update on the resulting reduced objective, for a total of \(400\) updates.
The IFeF fixed-coefficient baseline computes the lower-level response with the
GN-type solver and performs \(40\) updates using the partial derivative
\(\nabla_{\boldsymbol\omega}\mathcal J(\widehat{\boldsymbol\beta}^k,\boldsymbol\omega^k)\)
with the fitted coefficients held fixed. The covariance
matrix adaptation evolution strategy (CMA-ES) baseline, implemented in
Optuna~\cite{akiba2019optuna} with its default configuration, minimizes the
validation objective \(\mathcal J\) over \(40\) trials, with search ranges
\(\log R_m,\log\alpha_B\in[-10,10]\) and \(\log\lambda\in[-20,5]\).

\begin{figure}[!htb]
    \centering
    \subcaptionbox{Hypergradient error\label{fig:exampleB-hypergradient-error-panel}}[0.48\textwidth]{%
        \includegraphics[width=\linewidth]{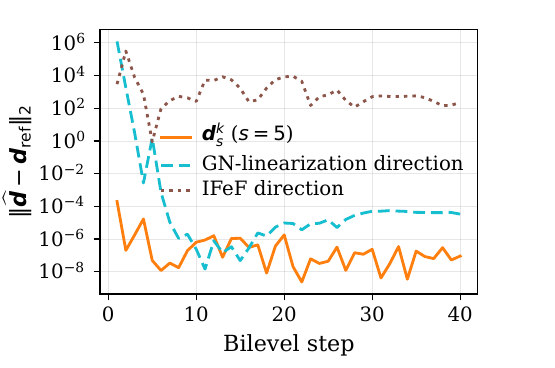}}
    \hfill
    \subcaptionbox{GN curvature diagnostic\label{fig:exampleB-gn-residual-panel}}[0.48\textwidth]{%
        \includegraphics[width=\linewidth]{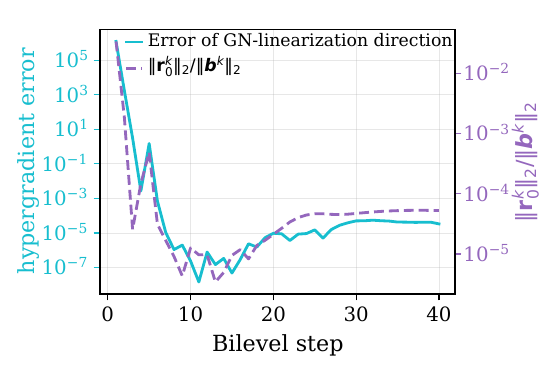}}

    \par\medskip
    \subcaptionbox{Best-so-far error vs.\ wall-clock time\label{fig:exampleB-wallclock-panel}}[0.62\textwidth]{%
        \includegraphics[width=\linewidth]{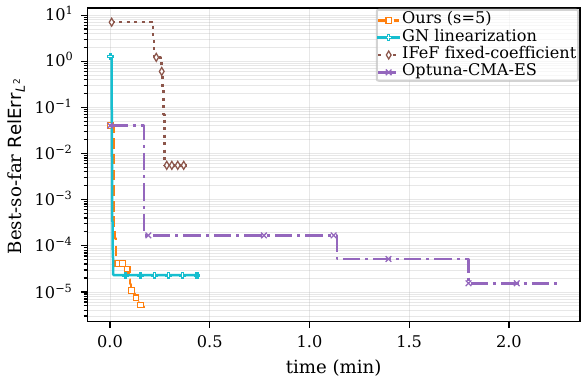}}
    \caption{Hypergradient errors versus bilevel step and best-so-far error versus wall-clock time
    for Example B with \(M=1000\).}
    \label{fig:exampleB-search-efficiency}
\end{figure}

We first examine the accuracy of the hypergradient directions for the three
derivative-based methods. As in Example~A, we compare \(\boldsymbol d_s^k\), the GN-linearization direction, and the IFeF
direction with exact-Hessian references at their respective computed responses
in Figure~\ref{fig:exampleB-search-efficiency}(a). The error of \(\boldsymbol d_s^k\) is much smaller than those
of the IFeF and GN-linearization directions. The large IFeF direction error results from omitting
the implicit dependence of the fitted coefficients on \(\boldsymbol\omega\).
In particular, \(\alpha_B\) and \(\lambda\) enter only the lower-level objective,
so the corresponding components of
\(\nabla_{\boldsymbol\omega}\mathcal J(\widehat{\boldsymbol\beta}^k,\boldsymbol\omega^k)\)
vanish and the IFeF baseline updates only \(R_m\).
To examine the effect of omitting the second-order term
\(\Delta_{\mathrm{GN}}^k\) in \eqref{eq:npde:local-hessian-decomposition},
Figure~\ref{fig:exampleB-search-efficiency}(b) compares the error of the
GN-linearization direction with the normalized initial adjoint residual
\(\|\mathbf r_0^k\|_2/\|\boldsymbol b^k\|_2
=\|\Delta_{\mathrm{GN}}^k\mathbf v_0^k\|_2/\|\boldsymbol b^k\|_2\) at each
GN-linearization response.
Their similar trends suggest that omitting \(\Delta_{\mathrm{GN}}^k\)
is an important source of this error.

Figure~\ref{fig:exampleB-search-efficiency}(c) and
Table~\ref{tab:exampleB-time-to-target} compare the four methods in terms of
solution accuracy and the time to reach each method's best error.
The proposed method attains the smallest error,
\(\operatorname{RelErr}_{L^2}=4.68\times10^{-6}\), in \(0.166\)~min. GN linearization reaches its best
error in \(0.017\)~min, but that error is about five times larger.
CMA-ES takes \(1.798\)~min to reach its best error of \(1.53\times10^{-5}\),
and IFeF yields the largest error.

\begin{table}[!htb]
\centering
\small
\setlength{\tabcolsep}{4pt}
\begin{tabular}{lcc}
\hline
Method
& Best-so-far \(\operatorname{RelErr}_{L^2}\)
& Time to best-so-far (min) \\
\hline
Ours, \(s=5\) & \(\mathbf{4.6848\times10^{-6}}\) & \(0.166\) \\
GN linearization & \(2.3066\times10^{-5}\) & \(0.017\) \\
IFeF fixed-coefficient & \(5.5375\times10^{-3}\) & \(0.272\) \\
Optuna--CMA-ES & \(1.5335\times10^{-5}\) & \(1.798\) \\
\hline
\end{tabular}
\caption{Best-so-far \(\operatorname{RelErr}_{L^2}\) and wall-clock time for Example B with \(M=1000\).}\label{tab:exampleB-time-to-target}
\end{table}

\begin{figure}[!htb]
    \centering
    \subcaptionbox{Exact solution\label{fig:exampleB-solution-exact}}[0.32\textwidth]{%
        \includegraphics[width=\linewidth]{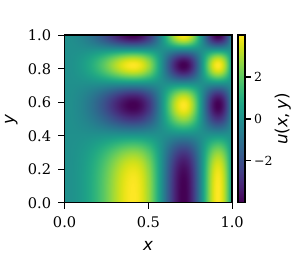}}
    \hfill
    \subcaptionbox{RFM/ELM solution\label{fig:exampleB-solution-elm}}[0.32\textwidth]{%
        \includegraphics[width=\linewidth]{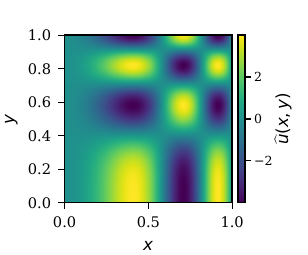}}
    \hfill
    \subcaptionbox{Absolute error\label{fig:exampleB-solution-error}}[0.32\textwidth]{%
        \includegraphics[width=\linewidth]{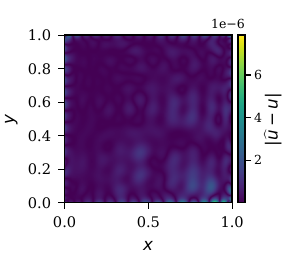}}
    \caption{Large-scale evaluation after \(100\)
    bilevel steps for Example B with \(M=2000\).}
    \label{fig:exampleB-solution}
\end{figure}

Finally, we test the proposed method in the larger benchmark setting of
SO-PIFRNN~\cite{linghu2025so}, with \(M=2000\) random features and \(2995\)
interior and \(2974\) boundary collocation points in each of the training and
validation sets. We initialize the solver-design parameters \(\boldsymbol\omega\)
at the values selected in the \(M=1000\) experiment. As shown in
Table~\ref{tab:exampleB-so-refit}, holding \(\boldsymbol\omega\) fixed and
refitting only the coefficients \(\boldsymbol\beta\) yields
\(\operatorname{RelErr}_{L^2}=1.29\times10^{-6}\), below the \(1.494\times10^{-6}\) reported for
SO-PIFRNN. Starting from the same \(\boldsymbol\omega\), \(100\) bilevel
steps reduce both \(\operatorname{RelErr}_{L^2}\) and
\(\operatorname{Err}_{L^\infty}\) by about a factor of five relative to this
coefficient-only refit. Figure~\ref{fig:exampleB-solution} shows the final solution
and its pointwise absolute error.

\begin{table}[!htb]
\centering
\small
\begin{tabular}{lccccc}
\hline
Evaluation state & \(R_m\) & \(\alpha_B\) & \(\lambda\) & $\operatorname{RelErr}_{L^2}$ & $\operatorname{Err}_{L^\infty}$ \\
\hline
Coefficient-only refit
& \(28.7269\) & \(65.2078\) & \(1.1548\times10^{-5}\)
& \(1.2861\times10^{-6}\) & \(4.2048\times10^{-5}\) \\
After \(100\) bilevel steps
& \(33.1457\) & \(78.4373\) & \(3.0922\times10^{-6}\)
& \(\mathbf{2.7282\times10^{-7}}\) & \(7.8777\times10^{-6}\) \\
\hline
\end{tabular}
\caption{Large-scale (\(M=2000\)) evaluation under the SO-PIFRNN sampling setting for Example B.}\label{tab:exampleB-so-refit}
\end{table}

\FloatBarrier 

\subsection{Example C (NeuroKansa): Burgers equation with adaptive centers}
\label{subsec:exampleC}

Examples~A and~B optimize global solver-design parameters
while keeping the center locations or sampled base parameters fixed. The third example learns a center-dependent Gaussian lengthscale over
space--time and uses it to guide center enrichment. We consider the viscous Burgers benchmark from
\cite{hu2026learning} on \((t,x)\in(0,4)\times(-10,10)\):
\[
u_t(t,x)+u(t,x)u_x(t,x)-\nu u_{xx}(t,x)=0,
\qquad \nu=0.5.
\]
The initial condition is \(u(0,x)=u_{\mathrm{ref}}(0,x)\), where
\(u_{\mathrm{ref}}(t,x)=\frac{1}{2}-\frac{1}{2}\tanh\!\left(\frac{x-\frac{1}2 t}{2}\right)\)
is a traveling-wave solution. At the truncated boundaries, we impose
\(u(t,-10)=1\) and \(u(t,10)=0\).
Since \(u_{\mathrm{ref}}\) does not exactly satisfy these boundary conditions,
we use it only as the reference solution for computing
\(\operatorname{RelErr}_{L^2}\) on the evaluation grid.

For a set of \(M\) space--time centers
\(\boldsymbol c_j=(t_j,x_j)\), the NeuroKansa approximation is
\[
\widehat{u}(t,x;\boldsymbol\beta,\boldsymbol\omega)
=
\sum_{j=1}^{M}\beta_j
\exp\!\left(
-\frac{\|(t,x)-\boldsymbol c_j\|_2^2}{\ell(\boldsymbol c_j;\boldsymbol\omega)^2}
\right).
\]
We
parameterize the center-dependent field \(\ell(\cdot;\boldsymbol\omega)\) by a two-hidden-layer multilayer
perceptron of width \(16\), with sigmoid linear unit activations and
Xavier-initialized hidden weights, whose scalar output is passed
through a positive map, so that \(\ell(\boldsymbol c;\boldsymbol\omega)>0\)
by construction.

Because the lengthscale field \(\ell(\cdot;\boldsymbol\omega)\) is learned
during bilevel optimization, we use its current values to guide the placement of
new interior centers at steps \(50\) and \(100\). At each of these steps \(k\), we generate a pool of \(408\) interior
candidates \(\{\boldsymbol z_i\}_{i=1}^{408}\) and evaluate
\(\ell_i^k=\ell(\boldsymbol z_i;\boldsymbol\omega^k)\).
With \(\ell_{\min}^k=\min_i\ell_i^k\) and
\(\ell_{\max}^k=\max_i\ell_i^k\), each candidate is accepted independently
with probability
\[
p_i^k=1-
\frac{\ell_i^k-\ell_{\min}^k}{\max\{\ell_{\max}^k-\ell_{\min}^k,\varepsilon_{\rm samp}\}},
\]
where \(\varepsilon_{\rm samp}=10^{-12}\) safeguards the denominator,
and the first \(51\) accepted candidates in pool order are appended as new
centers. Since \(p_i^k\) decreases with \(\ell_i^k\), new centers are
preferentially placed where the current representation uses more localized
basis functions.

\begin{figure}[!htb]
    \centering
    \subcaptionbox{Reference solution\label{fig:exampleC-reference}}[0.32\textwidth]{%
        \includegraphics[width=\linewidth]{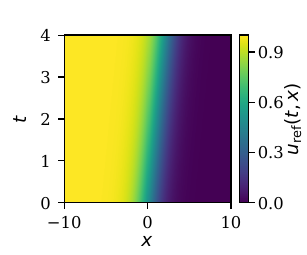}}
    \hfill
    \subcaptionbox{NeuroKansa solution\label{fig:exampleC-neurokansa}}[0.32\textwidth]{%
        \includegraphics[width=\linewidth]{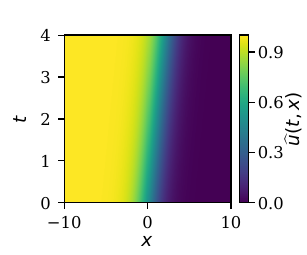}}
    \hfill
    \subcaptionbox{Absolute error\label{fig:exampleC-error}}[0.32\textwidth]{%
        \includegraphics[width=\linewidth]{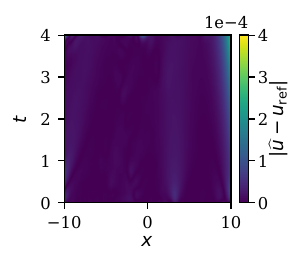}}
    \par\medskip
    \subcaptionbox{Bilevel trace\label{fig:exampleC-trace}}[0.32\textwidth]{%
        \includegraphics[width=\linewidth]{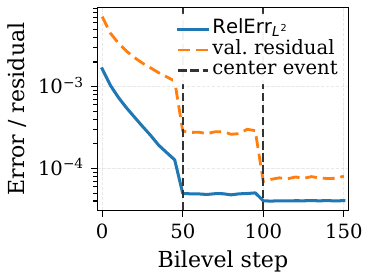}}
    \hfill
    \subcaptionbox{Reference slices\label{fig:exampleC-slices}}[0.32\textwidth]{%
        \includegraphics[width=\linewidth]{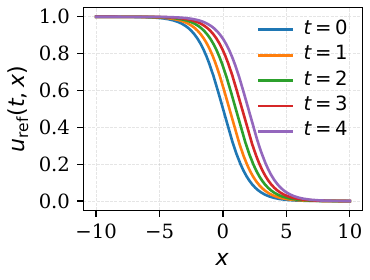}}
    \hfill
    \subcaptionbox{Learned lengthscale field\label{fig:exampleC-lengthscale}}[0.32\textwidth]{%
        \includegraphics[width=\linewidth]{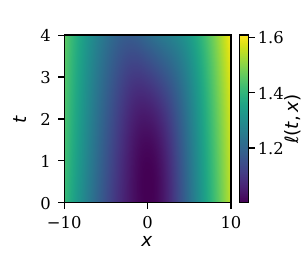}}
    \caption{Adaptive NeuroKansa solution and diagnostics after
    \(150\) bilevel steps for Example C. The vertical lines in panel~(d) mark the two scheduled center
    additions.}
    \label{fig:exampleC-summary}
\end{figure}

We use the bilevel formulation
\eqref{eq:outer_objective}--\eqref{eq:inner_constraint} with
\(\mathbf W_{\mathcal T}=\mathbf I\) and \(\lambda=0\).
The initial trial space
contains \(512\) uniformly sampled interior centers and \(96\) fixed
initial/boundary centers, with \(64\) initial centers and \(16\) on each
spatial boundary. The training and validation sets each contain \(1024\)
interior points with the same \(64/16/16\) initial/boundary allocation, and
errors are evaluated on a uniform \(256\times256\) space--time grid. We run
\(150\) bilevel steps with Adam, using initial step size \(0.002\) and
multiplicative decay \(0.99\). The adjoint computation uses \(\mu=10^{-9}\),
\(\tau_{\mathrm{adj}}=10^{-8}\), and \(s=5\).
Panels~(a)--(c) of Figure~\ref{fig:exampleC-summary} show the reference
solution, the NeuroKansa solution after \(150\) bilevel steps, and its
pointwise absolute error, respectively. In panel~(d), both the validation
residual and \(\operatorname{RelErr}_{L^2}\) decrease over the run, with drops
at the two center additions.

The reference profiles in Figure~\ref{fig:exampleC-summary}(e) show the
transition layer moving through the central part of the domain. In panel~(f),
the learned lengthscales are smaller in this region and larger toward the
spatial boundaries. Table~\ref{tab:exampleC-reference-lengthscale} compares the solution errors
obtained with the learned lengthscale field and the fixed scalar lengthscale
\(\ell=1/0.9\) used in \cite{hu2026learning}. For the fixed lengthscale, we
retain the final centers from the bilevel run and refit the coefficients
\(\boldsymbol\beta\) by solving the lower-level problem on the same training
set. The learned field yields
\(\operatorname{RelErr}_{L^2}=4.05\times10^{-5}\), about one fifth of the
value obtained with the fixed lengthscale (\(1.93\times10^{-4}\)).

\begin{table}[!htb]
\centering
\small
\begin{tabular}{lc}
\hline
Lengthscale choice & $\operatorname{RelErr}_{L^2}$ \\
\hline
Fixed scalar \(\ell=1/0.9\)~\cite{hu2026learning}
& \(1.9297\times10^{-4}\) \\
Learned \(\ell(\boldsymbol c;\boldsymbol\omega)\)
& \(\mathbf{4.0531\times10^{-5}}\) \\
\hline
\end{tabular}
\caption{Final \(\operatorname{RelErr}_{L^2}\) for the learned field and
the fixed scalar lengthscale, using the same
final center and collocation sets, for Example C.}\label{tab:exampleC-reference-lengthscale}
\end{table}

\FloatBarrier 

\subsection{Example D (IFeF-PINN): Bilevel feature refinement for Burgers equation}
\label{subsec:exampleD}

We consider the viscous Burgers benchmark from IFeF-PINN~\cite{wu2026iterative},
where earlier training stages provide a warm-start network state
for a final feature-refinement stage. Retaining this warm-start state, we
reformulate the final stage as a validation-driven bilevel problem. The
governing equation is
\[
    u_t + u u_x - \nu u_{xx}=0,\qquad
    (t,x)\in(0,0.99)\times(-1,1),\qquad
    \nu=\frac{0.01}{\pi}.
\]
The initial condition is \(u(0,x)=-\sin(\pi x)\), with homogeneous
Dirichlet boundary conditions \(u(t,-1)=u(t,1)=0\).

Let
\(h_{\boldsymbol\omega}(t,x)\in\mathbb R^{d_h}\) denote the hidden
representation. The fixed random Fourier matrix
\(W\in\mathbb R^{d_h\times N_F}\) has entries sampled independently from
\(\mathcal N(0,\sigma^2)\). Define
\[
\boldsymbol\psi_{\boldsymbol\omega}(t,x)
=
\sqrt{\frac{2}{N_F}}
\begin{bmatrix}
\cos\bigl(W^\top h_{\boldsymbol\omega}(t,x)\bigr)\\
\sin\bigl(W^\top h_{\boldsymbol\omega}(t,x)\bigr)
\end{bmatrix}
\in\mathbb R^{2N_F},
\]
and write
\[
\widehat{u}(t,x;\boldsymbol\beta,\boldsymbol\omega)
=
\boldsymbol\psi_{\boldsymbol\omega}(t,x)^\top\boldsymbol\beta.
\]
Here \(\boldsymbol\omega\) contains the trainable hidden-layer parameters,
whereas \(\boldsymbol\beta\in\mathbb R^{2N_F}\) is the coefficient vector of
the linear output layer. We take \(d_h=20\), \(N_F=800\), and \(\sigma=1\). The
remaining network settings are the same as in \cite{wu2026iterative}.

We use the bilevel formulation
\eqref{eq:outer_objective}--\eqref{eq:inner_constraint} with
\(\mathbf W_{\mathcal T}=\mathbf I\), \(\mathbf R_{\mathcal T}=\mathbf I\), and
\(\lambda=2\times10^{-7}\). Following the relative residual weighting in the Burgers implementation of
\cite{wu2026iterative}, we scale the interior and initial/boundary blocks of
\(\boldsymbol\rho(\boldsymbol\beta,\boldsymbol\omega;\mathcal C)\) by
\(\sqrt{0.02}\) and \(\sqrt{2}\), respectively, for each
\(\mathcal C\in\{\mathcal T,\mathcal V\}\).

\begin{figure}[!htb]
    \centering
    \subcaptionbox{Reference solution\label{fig:exampleD-ours-exact}}[0.32\textwidth]{%
        \includegraphics[width=\linewidth]{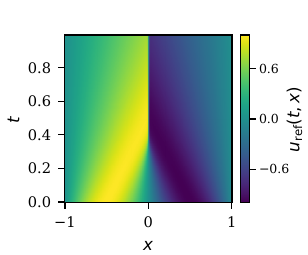}}
    \hfill
    \subcaptionbox{Bilevel solution\label{fig:exampleD-ours-prediction}}[0.32\textwidth]{%
        \includegraphics[width=\linewidth]{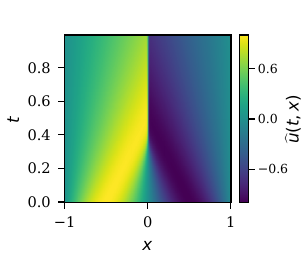}}
    \hfill
    \subcaptionbox{Absolute error\label{fig:exampleD-ours-error}}[0.32\textwidth]{%
        \includegraphics[width=\linewidth]{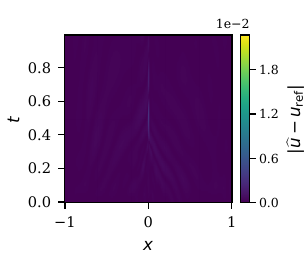}}
    \caption{Bilevel solution after \(500\) bilevel steps for Example D.}
    \label{fig:exampleD-ours-solution}
\end{figure}

The training and validation sets are sampled independently and each contain \(10{,}000\)
interior and \(100\) initial/boundary points.
Errors are evaluated against the numerical reference solution provided with
\cite{raissi2019physics} on its uniform \(100\times256\) grid in \((t,x)\).
Starting from the warm-start state, we run \(500\) bilevel steps with AdamW, using step size
\(10^{-4}\) and weight decay \(10^{-5}\). The
adjoint computation uses
\(\mu=10^{-6}\),
\(\tau_{\mathrm{adj}}=10^{-6}\), and \(s=5\). The experiment was executed on an
NVIDIA A100 GPU.
Figure~\ref{fig:exampleD-ours-solution} shows the
reference solution, the bilevel solution after \(500\) bilevel steps, and
its pointwise absolute error.

\begin{table}[!htb]
\centering
\small
\begin{tabular}{lcc}
\hline
Method & Final-stage updates & $\operatorname{RelErr}_{L^2}$ \\
\hline
IFeF-PINN baseline & \(2000\) & \(3.1005\times10^{-3}\) \\
Ours, \(s=5\) & \(500\) & \(\mathbf{3.0983\times10^{-4}}\) \\
\hline
\end{tabular}
\caption{Final-stage \(\operatorname{RelErr}_{L^2}\) for IFeF-PINN and the proposed method for Example D.}\label{tab:exampleD-ifef}
\end{table}

\begin{figure}[!htb]
    \centering
    \subcaptionbox{\(\operatorname{RelErr}_{L^2}\) trajectory\label{fig:exampleD-ifef-trajectory}}[0.32\textwidth]{%
        \includegraphics[width=\linewidth]{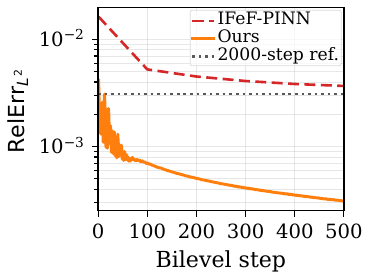}}
    \hfill
    \subcaptionbox{IFeF-PINN solution\label{fig:exampleD-ifef-solution}}[0.32\textwidth]{%
        \includegraphics[width=\linewidth]{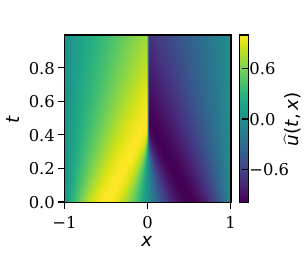}}
    \hfill
    \subcaptionbox{Absolute error\label{fig:exampleD-ifef-error}}[0.32\textwidth]{%
        \includegraphics[width=\linewidth]{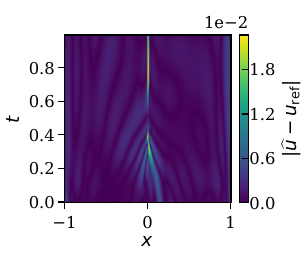}}
    \caption{Final-stage comparison with IFeF-PINN for Example D.
    The horizontal dotted line in panel~(a) marks the IFeF-PINN endpoint \(\operatorname{RelErr}_{L^2}\)
    after \(2000\) updates.}
    \label{fig:exampleD-ifef}
\end{figure}

To evaluate the proposed method, we compare it with the final-stage refinement
of the IFeF-PINN baseline~\cite{wu2026iterative}. IFeF-PINN performs \(2000\) final-stage
updates of \(\boldsymbol\omega\), each using the gradient of the training loss with
\(\boldsymbol\beta\) held fixed, and refits \(\boldsymbol\beta\) every \(100\)
updates. As shown in Table~\ref{tab:exampleD-ifef} and
Figure~\ref{fig:exampleD-ifef}(a), the proposed method attains
\(\operatorname{RelErr}_{L^2}=3.10\times10^{-4}\) after \(500\) bilevel steps,
an order of magnitude smaller than the \(3.10\times10^{-3}\) attained by
IFeF-PINN after \(2000\) updates. Comparing
Figures~\ref{fig:exampleD-ours-solution}(b)--(c) and
\ref{fig:exampleD-ifef}(b)--(c) shows that the proposed method yields smaller
pointwise errors along the steep transition at \(x=0\), particularly at later
times.

\FloatBarrier 
\subsection{Example E (RFM/ELM): Five-dimensional nonlinear Helmholtz problem}
\label{subsec:exampleE}

To test the proposed method on a higher-dimensional problem, we extend the nonlinear
Helmholtz problem \eqref{eq:helmholtz} in Example~B to \(\Omega=(0,1)^5\).
The Dirichlet boundary data are prescribed
by \(g(\boldsymbol x)=u(\boldsymbol x)\) for
\(\boldsymbol x\in\partial\Omega\), where the exact solution is
\[
u(\boldsymbol x)=4\prod_{i=1}^{5}\cos\!\left(\frac{\pi x_i}{2}\right),
\]
and the forcing term is chosen accordingly as
\(
f(\boldsymbol x)=\Delta u(\boldsymbol x)-100u(\boldsymbol x)
+10\cos\bigl(u(\boldsymbol x)\bigr).
\)
We use the sine-feature approximation and bilevel formulation of
Example~B, with \(\widetilde{\boldsymbol w}_j\sim\mathcal U([-1,1]^5)\) and
\(\boldsymbol\omega=(\log R_m,\log\alpha_B,\log\lambda)\).

\begin{figure}[!htb]
    \centering
    \subcaptionbox{Exact solution\label{fig:exampleE-solution-exact}}[0.32\textwidth]{%
        \includegraphics[width=\linewidth]{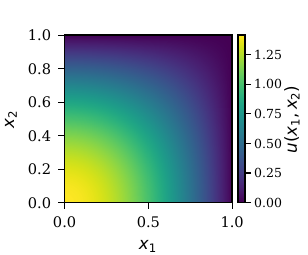}}
    \hfill
    \subcaptionbox{RFM/ELM solution\label{fig:exampleE-solution-elm}}[0.32\textwidth]{%
        \includegraphics[width=\linewidth]{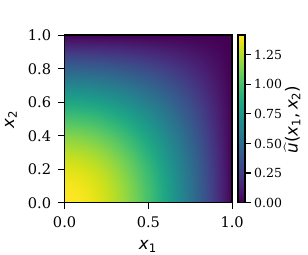}}
    \hfill
    \subcaptionbox{Absolute error\label{fig:exampleE-solution-error}}[0.32\textwidth]{%
        \includegraphics[width=\linewidth]{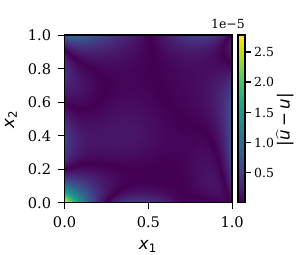}}
    \caption{RFM/ELM solution after \(80\)
    bilevel steps on the slice \(x_3=x_4=x_5=\tfrac12\) for Example E.}
    \label{fig:exampleE-solution}
\end{figure}

The experiment uses \(M=2000\) sine random features.
The training and validation sets each contain
\(1000\) interior and \(1000\) boundary points, and errors are
evaluated on an independently sampled set of \(7000\) interior and \(1000\)
boundary points. Starting from \((R_m,\alpha_B,\lambda)=(1,1,10^{-8})\), we
run \(80\) bilevel steps with Adam, using constant step size \(0.1\).
The adjoint computation uses
\(\mu=10^{-10}\), \(\tau_{\mathrm{adj}}=10^{-8}\), and \(s=5\).
Figure~\ref{fig:exampleE-solution} shows the exact solution, the
RFM/ELM solution after \(80\) bilevel steps, and the
pointwise absolute error on the two-dimensional slice
\(x_3=x_4=x_5=\tfrac12\).

To evaluate the proposed method, we
compare it with two baselines implemented in
Optuna~\cite{akiba2019optuna}: the tree-structured Parzen estimator (TPE) and
CMA-ES. Both
use the same initial parameters, training/validation sets, and number of random features, and minimize
the validation objective \(\mathcal J\) over \(80\) trials, with search ranges \(\log R_m,\log\alpha_B\in[-10,10]\)
and \(\log\lambda\in[-20,5]\).

\begin{figure}[!htb]
    \centering
    \includegraphics[width=0.62\textwidth]{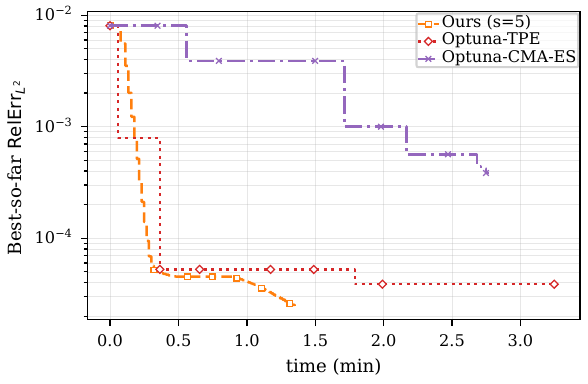}
    \caption{Best-so-far error versus wall-clock time for Example E with \(M=2000\).}
    \label{fig:exampleE-wallclock}
\end{figure}

We assess the computational efficiency and solution accuracy
of the three methods in Figure~\ref{fig:exampleE-wallclock} and
Table~\ref{tab:exampleE-search}. The proposed method attains the smallest error with the
shortest search time, whereas TPE requires more than twice the search time
and yields a larger error, and CMA-ES is about an order of magnitude
less accurate.

\begin{table}[!htb]
\centering
\small
\setlength{\tabcolsep}{4pt}
\begin{tabular}{lccc}
\hline
Method & Best-so-far \(\operatorname{RelErr}_{L^2}\)
& Time to best-so-far (min) & Search time (min) \\
\hline
Ours, \(s=5\) & \(\mathbf{2.4817\times10^{-5}}\) & \(1.349\) & \(1.349\) \\
Optuna--TPE & \(3.8793\times10^{-5}\) & \(1.789\) & \(3.272\) \\
Optuna--CMA-ES & \(3.8436\times10^{-4}\) & \(2.748\) & \(2.770\) \\
\hline
\end{tabular}
\caption{Best-so-far \(\operatorname{RelErr}_{L^2}\) and
wall-clock time for Example E with \(M=2000\).}
\label{tab:exampleE-search}
\end{table}

\FloatBarrier

\section{Conclusions and Perspectives}
\label{sec:conclusions}

We consider several classes of mesh-free partial differential
equation (PDE) solvers whose formulations share a common structure:
coefficients are fitted by solving a least-squares problem, while solver
accuracy depends strongly on the choice of solver-design parameters. Building
on this shared structure, we introduce a unified validation-driven bilevel optimization
method that separates training-based coefficient fitting from validation-based
solver-design selection. Implicit hypergradient methods optimize the reduced validation
objective along a selected branch of lower-level stationary points.
For nonlinear PDEs, evaluating these hypergradients is
computationally demanding. To make this approach efficient,
we develop an approximate implicit hypergradient method that combines
finite warm-started damped Gauss--Newton (GN) lower-level solves with
GN-preconditioned minimum residual iterations for the adjoint system. To
characterize how these computational approximations affect the bilevel
optimization and ultimately the accuracy of the resulting PDE solver, we
establish an error analysis connecting hypergradient approximation and
validation sampling to PDE solution accuracy. Numerical experiments on nonlinear
PDEs demonstrate the applicability of the framework across Kansa methods,
random feature and extreme learning machine methods, and methods based on
physics-informed neural networks. The numerical comparisons also show that the proposed
method achieves accurate solutions with favorable time-to-accuracy performance.

Taken together, these results establish a unified framework for
validation-driven design of mesh-free PDE solvers, linking
practical bilevel optimization with quantitative control of PDE solution
accuracy. The framework also suggests several extensions toward broader PDE
classes and more complex learning
problems.
\begin{itemize}
    \item The optimizer-independent terminal PDE error bound does not require a
    Polyak--\L{}ojasiewicz condition, but the finite-iteration analysis of the
    basic update does. Establishing trajectory bounds under weaker conditions
    on the upper-level objective remains an important direction.

    \item The current validation-to-PDE result relies on uniform
    bounds on validation sampling error, norm compatibility, and residual stability. The next step is to
    investigate how broadly this connection holds across different PDE classes,
    solution norms, and validation settings.

    \item The current implementations and analysis focus on lower-level
    responses defined by coefficient fitting. In more general PDE-constrained
    optimization and inverse-problem settings, the lower level may instead
    involve training a
    full neural network \cite{haobi,zhang2026bilo}. In such settings, the
    lower-level optimization is substantially more complex and less stable,
    making it difficult to converge rapidly to the accuracy required for a
    reliable hypergradient. A central question is therefore how to balance
    lower-level solve accuracy and computational cost while controlling the resulting
    hypergradient error.
\end{itemize}

\apptocmd{\thebibliography}{%
    \setlength{\itemsep}{-1pt}%
}{}{}
{\small
\bibliographystyle{abbrv}
\bibliography{references.bib}
}
\end{document}